\RequirePackage{fix-cm}
\documentclass[smallextended, envcountsect]{svjour3}       
\smartqed  
\usepackage{graphicx}
\usepackage{latexsym,amsmath,amssymb,amsfonts}
\usepackage[mathscr]{eucal}
\usepackage{color}
\usepackage{soul,xcolor}

\newtheorem{theo}{Theorem}[section]
\newtheorem{cor}{Corollary}[section]

\newcommand{\cd}{\mathcal{D}}
\newcommand{\cl}{\mathcal{L}}
\newcommand{\cn}{\mathcal{N}}
\newcommand{\ca}{\mathcal{A}}
\newcommand{\cu}{\mathcal{U}}
\newcommand{\ch}{\mathcal{H}}
\newcommand{\cg}{\mathcal{G}}
\newcommand{\ck}{\mathcal{K}}
\newcommand{\cm}{\mathcal{M}}
\newcommand{\mmin}{\mu_{min}}
\newcommand{\mmax}{\mu_{max}}

\definecolor{ForestGreen}{cmyk}{0.91,0,0.88,0.12}

\newcommand{\epsorpdf}[1]{eps/#1.eps} 

\begin{document}

\title{A comparison of the Extrapolated Successive Overrelaxation and the Preconditioned Simultaneous Displacement methods for augmented linear systems
}

\titlerunning{A comparison of the ESOR and PSD methods for augmented linear systems}        

\author{M. A. Louka \and N. M. Missirlis}

\authorrunning{M. A. Louka et al.} 

\institute{M. A. Louka \and N. M. Missirlis \at
              Department of Informatics and Telecommunications, University of Athens, \\
        Panepistimiopolis, 15784, Athens, Greece\\
              Tel.: +30-72-75103\\
              Fax: +30-72-75114\\
              \email{mlouka@di.uoa.gr; nmis@di.uoa.gr }           
}

\date{Received: date / Accepted: date}

\maketitle

\begin{abstract}
In this paper we study the impact of two types of preconditioning on the numerical solution of large sparse augmented linear systems. The first preconditioning matrix is the lower triangular part whereas the second is the product of the lower triangular part with the upper triangular part of the augmented system's coefficient matrix. For the first preconditioning matrix we form the Generalized Modified Extrapolated Successive Overrelaxation (GMESOR) method, whereas the second preconditioning matrix yields the Generalized Modified Preconditioned Simultaneous Displacement (GMPSD) method, which is an extrapolated form of the Symmetric Successive Overrelaxation method. We find sufficient conditions for each aforementioned
iterative method to converge. In addition, we develop a geometric approach, for determining the optimum values of their parameters and corresponding spectral radii. It is shown that both iterative methods studied (GMESOR and GMPSD) attain the same rate of convergence. Numerical results confirm our theoretical expectations.
\keywords{Iterative methods \and linear systems  \and augmented systems \and SOR.}
\subclass{AMS(MOS) \and 65F10 \and 65N20 \and CR:5.13.}
\end{abstract}

\section{Introduction}
\label{intro}
Let $A \in \mathbb{R}^{m\times m}$ be a symmetric positive definite matrix and $B\in\mathbb{R}^{m\times n}$ be a matrix of full column rank, where $m\geq n$.  Then, the augmented linear system is of the form 

\begin{equation}
\ca u=b\label{eq:01}
\end{equation}
where
\begin{equation}\label{eq:02}
\ca=\left(
\begin{array}{cc}
 A & B\\
-B^{T} & 0
\end{array}
\right), \;\; u=
\left(
\begin{array}{cc}
x\\
y
\end{array}
\right),\;\; b=
\left(
\begin{array}{cc}
b_1\\
-b_2
\end{array}
\right)
\end{equation}
with $B^T$ denoting the transpose of the matrix $B$.
\\
When $A$ and $B$ are large and sparse matrices, iterative methods for solving  (\ref{eq:01})-(\ref{eq:02}) are effective and more attractive than direct methods, because of storage requirements and preservation of sparsity.
There are several approaches to the iterative solution of (\ref{eq:01})-(\ref{eq:02}). First, we mention multigrid methods \cite{BenGolLie2005}, {{\cite{Metsch2013},}} which are often the most efficient and effective methods for solving large, sparse, linear systems \cite{BriHe2000}, \cite{TrotOosSch2000}.
For example, one can apply multigrid techniques to the whole system (\ref{eq:01})-(\ref{eq:02}) to solve problems in areas of computational fluid dynamics \cite{GaLiOoVa07}, \cite{HaBeHa2010}, \cite{OOstGasp2008}, \cite{Vo2006}, \cite{WangChen2013}, \cite{GR98}, \cite{GrOeSch03} constrained optimization \cite{ScWi08}, \cite{SiZu09}, \cite{TaZu2009}, \cite{TaZu11}, mixed finite elements \cite{ArGo2009}, \cite{GoPo2012} and elsewhere. For parallel multigrid see e.g \cite{GrMeOeSc05}, \cite{GrMeSch06}, \cite{GrMeSch08}.


On the other hand the difficulty in applying iterative methods such as the Successive Overrelaxation (SOR) method \cite{Young} to the system (\ref{eq:01})-(\ref{eq:02}) is the singularity of the block diagonal part of the coefficient matrix. Various methods have been developed to overcome this problem such as the Uzawa and the Preconditioned Uzawa methods \cite{ArHurUz58}, \cite{BraPasVas77}, \cite{ElmGol94}. In 2001, Golub et al. \cite{GolWuYuan2001} generalized the Uzawa and the Preconditioned Uzawa methods by introducing an additional acceleration parameter and produced the SOR-like method. When a good preconditioning matrix is easily computed one can consider the MINRES and GMRES methods \cite{FiRaSiWa98}, \cite{BraPasVas77} for solving (\ref{eq:01})-(\ref{eq:02}). In case the matrix $A$ in (\ref{eq:02}) is symmetric and positive definite, the Preconditioned Conjugate Gradient (PCG) method \cite{HestSt52} can be applied. This was done with an SOR-like preconditioner in the work by Li, Evans and Zhang in \cite{LiEvansZha2003}. In 2005, Bai et al. \cite{BaiParlWang2005} studied the Generalized SOR (GSOR) method by introducing an additional parameter to the SOR-like method and proved that it possesses the same rate of convergence but lower complexity than the PCG method.
Furthermore, the Generalized Modified Extrapolated SOR (GMESOR) method was also proposed for further study. The latter is a generalization of the GSOR method as it uses one additional parameter. The way of introducing parameters resembles the one followed for the formulation of the Modified SOR method \cite{Young}, \cite{Hadj}, \cite{MisEv1}, \cite{Mis1}, \cite{MisEv2} in case  of two-cyclic linear systems.
\\
\noindent The present paper was motivated by the work in \cite{BaiParlWang2005}. {{We develop the convergence analysis of the Generalized Modified Extrapolated SOR (GMESOR) method and the Generalized Modified Preconditioned Simultaneous Displacement (GMPSD) method}}.
These methods introduce more parameters with the hope to further increase their rate of convergence.
The goal of our work was to study the impact of two different preconditioning matrices to the convergence rate of the associated iterative method for solving the augmented linear system (\ref{eq:01})-(\ref{eq:02}). First, we use the preconditioning matrix which is formed by the lower triangular part of $\ca$ and formulate the GMESOR method which is an extrapolated form of the GSOR method. Secondly, we consider as preconditioning matrix the product of the lower with the upper triangular part of $\ca$  and construct the GMPSD method. The reason for studying the latter form of preconditioning matrix is to obtain a better approximation to the matrix $\ca$ than the former and as such it is hoped to produce an iterative method with a faster rate of convergence. The construction of both methods resembles the one followed for the MESOR and MPSD methods studied in \cite{Mis1} and \cite{MisEv2}, respectively, for two-cyclic matrices.
Our starting point, for studying these iterative methods, is the derivation of functional relationships which relate the eigenvalues of their iteration matrices with those of the matrix $J=Q^{-1}B^TA^{-1}B$. Assuming that the matrix Q is symmetric positive or negative definite, the eigenvalues of the matrix $J$ are real and either positive or negative, respectively. Under these assumptions we find sufficient conditions for the convergence of the GMESOR and GMPSD methods and determine the optimum values of their parameters.
The study of GMESOR and GMPSD  becomes interesting as these methods can be used either as preconditioners to {{Krylov subspace methods \cite{BaiGolPan2004}, \cite{GolWuYuan2001}, \cite{HestSt52}, \cite{FiRaSiWa98}}} or as smoothers to multilevel methods \cite{BaMo2012}, \cite{BaShu2013}, \cite{Ba2013}. Traditionally, multigrid methods utilize stationary iterative methods (such as Jacobi, Gauss-Seidel ) to smooth out high-frequency errors and accelerate the convergence. In \cite{CaoJiMiss13} a semi-iterative method, namely the Chebyshev-Jacobi method, was used as smoother. Similarly, the GMESOR method or the GMPSD method in combination with semi-iterative techniques can be used as smoothers to accelerate the rate of convergence of multigrid methods. Recent work for the application of algebraic multigrids for saddle point systems is presented in \cite{Metsch2013} and the references therein.

The contributions of our paper can be summarized as follows.\\
  (i) We present a simple and unified approach for developing the convergence analysis of
   the GMESOR and GMPSD methods. In particular, we develop a geometrical approach for the determination of the optimum values of the parameters in GMESOR and GMPSD methods which is similar to Varga \cite{Var1} p. 111, for finding the optimum value of the parameter $\omega$ in SOR.  The difference, in our case, is that now the functional relationship contains more than two parameters and consequently we had to extend the proof of \cite{Var1}.  There is a different algebraic approach in  \cite{Young} pp. 279 for the determination of the two optimum values for $\omega$ and $\omega^{'}$ in the Modified SOR (MSOR) method which, with some additional modifications, will solve the problem as far as the GSOR method is concerned. Nevertheless, it is doubtful whether this approach works also for the determination of the optimum value for more than two parameters as is the case for the GMESOR and GMPSD methods. This is also the case if one adopts the  approach of  \cite{BaiParlWang2005}.\\
(ii) From our theoretical and experimental analysis it is shown that both aforementioned forms of preconditioning matrices have the same impact on the convergence rate of the induced iterative method for the numerical solution of the augmented linear system (\ref{eq:01})-(\ref{eq:02}). More specifically, the GMESOR and GMPSD methods attain the same convergence rate since their spectral radii are identical for the optimum values of their parameters. In particular, we show that GMESOR degenerates to GSOR, whereas a simplified version of the GMPSD method is identical to a backward form of the GSOR method.
\noindent
Furthermore, we compare the effectiveness of our methods in relation to the PHSS  \cite{BaiGolNg2003}, \cite{BaiGolPan2004}, \cite{BaiGolLi2005}, \cite{BaiGol2007}, \cite{BaiGolLi2007}, \cite{WangBai2001} {{and Krylov subspace methods \cite{SaSch86}, \cite{Saad03}, \cite{HVdV03}}}.
\\
\noindent The paper is organized as follows. In section 2 we study the convergence of the GMESOR method. In particular, we find sufficient conditions  for GMESOR to converge under the assumption that the eigenvalues of the $J$ matrix are real. We also determine optimum values for its parameters. A similar convergence analysis for the GMPSD method is developed in section 3. In section 4, we present our numerical results and finally in section 5 we state our remarks and conclusions.
\section{The Generalized Modified Extrapolated SOR (GMESOR) method}
\label{sec:1}
%
Let the coefficient matrix $\ca$ of (\ref{eq:01}) be defined by the splitting
\begin{equation}\label{eq:03}
\ca=\cd-\cl-\cu
\end{equation}
where
\begin{equation}\cd=
\left(
\begin{array}{ccc}
A && 0    \\
0 && Q \\
\end{array}\right),\
\cl=\left(
\begin{array}{ccc}
0 & &0    \\
B^T & &aQ \\
\end{array}\right),\
\cu=\left(
\begin{array}{ccc}
0 & &-B    \\
0 & &(1-a)Q \\
\end{array}\right),
\label{eq:04}
\end{equation}
with $Q\in\mathbb{R}^{n \times n}$ be a prescribed nonsingular and
symmetric matrix and ${a\in\mathbb{R}}$. Furthermore, we denote by $T$, the diagonal matrix $T=diag(\tau_1 I_m, \tau_2 I_n)$ with $\tau_1,\tau_2\in \mathbb{R}-{\{0\}},\; I_m \in\mathbb{R}^{m \times m}$ and $I_n\in\mathbb{R}^{n \times n}$ be identity matrices.
\\ \noindent For the numerical solution of (\ref{eq:01}), we consider the following iterative scheme
\begin{equation}\left(
\begin{array}{c}
x^{(k+1)}    \\
y^{(k+1)}\\
\end{array}\right)=\ch{(\tau_1, \tau_2)}
\left(
\begin{array}{c}
x^{(k)}    \\
y^{(k)}\\
\end{array}\right)+\eta{(\tau_1, \tau_2)}
\left(\begin{array}{c}
b_1    \\
-b_2\\
\end{array}\right)
\label{eq:05}
\end{equation}
where
\begin{equation}\label{eq:06}
\ch{(\tau_1, \tau_2)}=I-R^{-1}T\ca, \;\;\eta{(\tau_1, \tau_2)}=R^{-1}Tb,
\end{equation}
$R$ is a nonsingular matrix to be defined and $I=diag(I_m, I_n)$.\\ \noindent
In the sequel we consider two different types of preconditioning matrices $R$ and study the corresponding iterative methods derived by (\ref{eq:05}) and (\ref{eq:06}).

\subsection{The functional relationship}
\label{sec:2}
As a first step we consider the preconditioning matrix which is formed by the parameterized diagonal and lower triangular part of $\ca$
\begin{equation}\label{eq:06bb}
R=\cd-\Omega\cl,
\end{equation}
where $\Omega=diag(\omega_1 I_m, \omega_2 I_n)$ with $\omega_1,\omega_2\in \mathbb{R}$.
If $R$ is given by
(\ref{eq:06bb}), then (\ref{eq:06}) becomes

\begin{equation*}
\ch{(\tau_1,\tau_2,\omega_2, a)}=I-{(\cd-\Omega\cl)}^{-1}T\ca
\end{equation*}
or
\begin{equation}\label{eq:08}
\ch{(\tau_1,\tau_2,\omega_2, a)}={(\cd-\Omega\cl)}^{-1}[(I-T)\cd+(T-\Omega)\cl+T\cu]
\end{equation}
and
\begin{equation}\label{eq:09}
\eta{(\tau_1,\tau_2,\omega_2, a)}={(\cd-\Omega\cl)}^{-1}Tb.
\end{equation}
Note that the parameter $\omega_1$ is absent in $\ch$ and $\eta$. This is because the first $m$ rows of $\cl$ are zeros a fact which is carried over in matrix $\Omega\cl$ also. \\ The iterative scheme given by (\ref{eq:05}),(\ref{eq:06}),(\ref{eq:08}) and (\ref{eq:09})   will be referred to as the Generalized Modified Extrapolated SOR (GMESOR) method. In case $a=0$ this method was introduced in \cite{BaiParlWang2005} and proposed for further study.
In the sequel to distinguish the dependence of GMESOR upon the parameter $a$ we use the notation GMESOR($a$).\\
For $(\cd-\Omega\cl)^{-1}$ to exist we require
\begin{equation*}
\det(\cd-\Omega\cl)\neq 0.
\end{equation*}
Because of (\ref{eq:04})
\begin{equation*}
R=\cd-\Omega\cl=\left(
\begin{array}{ccccccc}
A & 0\\
-\omega_2 B^T & (1-a\omega_2)Q
\end{array}
\right).
\end{equation*}
Therefore,
\begin{equation*}
\det(\cd-\Omega\cl)=(1-a\omega_2)^{n}\det{(A)}\det{(Q)}\neq 0
\end{equation*}
or
\begin{equation}\label{eq:10b}
a\omega_2\neq 1
\end{equation}
since the matrix $A$ is symmetric positive definite and the matrix $Q$ is nonsingular. In the sequel we require (\ref{eq:10b}) to hold.\\ \noindent
The GMESOR($a$) method has the following algorithmic form.
\\\\ \noindent {\sc The GMESOR($a$) Method:}{\em{ Let $Q\in\mathbb{R}^{n\times n}$ be a nonsingular and symmetric matrix.
Given initial vectors $x^{(0)}\in\mathbb{R}^m$ and
$y^{(0)}\in\mathbb{R}^n$, and the parameters $\tau_1,\tau_2\neq 0,  \
\omega_2, a\in\mathbb{R}$ with $a\omega_2\neq 1$. For $k=0,1,2,...$ until the
iteration sequence
$\{({x^{(k)}}^T, {y^{(k)}}^T)^T\}$ is convergent, compute \\ \newline
\hspace{0.3cm}${x^{(k+1)}=(1-\tau_1)x^{(k)}+\tau_1 A^{-1}(b_1-By^{(k)})}, $\\
\hspace{0.5cm}$    {y^{(k+1)}=y^{(k)}+
\displaystyle\frac{1}{1-a\omega_2}Q^{-1}\left\{B^T[\omega_2x^{(k+1)}+(\tau_2-\omega_2)x^{(k)}]-\tau_2b_2\right\}},
$ \\ \newline
\noindent where Q is an approximation of the Schur complement matrix $B^TA^{-1}B$.}}
\\\\
\noindent For special values of its parameters GMESOR($a$) degenerates into known methods. Indeed, if $\omega=\tau_1=\tau_2=\omega_2$ and $a=0$ then GMESOR becomes the SOR-like method \cite{GolWuYuan2001}; if $\omega=\tau_1=\tau_2=\omega_2=1$ and $a=0$ then it becomes the preconditioned Uzawa method \cite{ElmGol94}; 
and if $\tau_1=\omega_1, \; \tau_2=\omega_2$ and $a=0$, then it becomes the GSOR method \cite{BaiParlWang2005}.
By comparing the algorithmic structures of GMESOR($a$) and {{GSOR}}, one can verify that the former has an additional matrix times a vector computation. Finally, if
\begin{equation}\label{eq:bai}
\tau_1=\omega, \ \frac{\omega_2}{1-a\omega_2}=\gamma \ \ \mbox{and}\ \  \frac{\tau_2}{1-a\omega_2}=\tau
\end{equation}
then the GMESOR($a$) method becomes the Generalized Inexact Accelerated Overrelaxation (GIAOR) method \cite{BaiParlWang2005} and if
\begin{equation}\label{eq:bai2}
\tau_1=\omega, \tau_2=\omega_2 \ \mbox{and}\ \  \frac{\tau_2}{1-a\omega_2}=\tau
\end{equation}
the {{GMESOR($a$)}} method becomes the Parametrized Inexact Uzawa (PIU) method \cite{BaiWang2008} when $P=A$.
\noindent The following theorem establishes the functional relationship between the eigenvalues $\lambda$ of the iteration matrix $\ch(\tau_1, \tau_2, \omega_2, a)$  and the eigenvalues $\mu$ of the associated matrix  $J=Q^{-1}B^TA^{-1}B$.
\begin{theo}\label{theo:func_rel}
Let $A\in\mathbb{R}^{m\times m}$ be symmetric positive definite, $B\in\mathbb{R}^{m\times n}$ be of full column rank and $Q\in\mathbb{R}^{n\times n}$ be
nonsingular and symmetric. If $\lambda\neq 1-\tau_1$ is an eigenvalue of the matrix
$\ch{(\tau_1,\tau_2,\omega_2,a)}$ of the GMESOR($a$) method and if $\mu$ satisfies
\begin{equation}
\lambda^2+\lambda\left(\tau_1-2+\frac{\tau_1\omega_2}{1-a\omega_2}\mu\right)+1-\tau_1+\frac{\tau_1(\tau_2-\omega_2)}{1-a\omega_2}\mu=0, \label{eq:10}
\end{equation}
where $a{\omega_2}\neq {1}$, then $\mu$ is an eigenvalue of the matrix $J=Q^{-1}B^TA^{-1}B$.
Conversely, if $\mu$ is an eigenvalue of J and if
$\lambda\neq1-\tau_1$ satisfies (\ref{eq:10}), then $\lambda$ is an
eigenvalue of $\ch{(\tau_1,\tau_2,\omega_2,a)}$. In addition, $\lambda=1-\tau_1$ is an
eigenvalue of $\ch{(\tau_1,\tau_2,\omega_2,a)}$ (if $m>n$) with the corresponding
eigenvector $(x^T,0)^T$, where $x\in\cn(B^T)$ and $\cn(B^T)$ is the null space of $B^T$.
\end{theo}
{ {{\bf{Proof}} The eigenvalues $\mu$ of the matrix $J=Q^{-1}B^TA^{-1}B$ are real, positive and non-zero. Let $\lambda$ be a nonzero eigenvalue of the iteration matrix $\ch(\tau_1, \tau_2, \omega_2,a)$ defined in (\ref{eq:08}), and $(x,y)^T\in\mathbb{R}^{m+n}$ be the corresponding eigenvector. Then,
\begin{equation}\label{eq:11a}
\ch(\tau_1, \tau_2, \omega_2, a)\left(
\begin{array}{c}
x\\y
\end{array}
\right)=\lambda\left(
\begin{array}{c}
x\\y
\end{array}
\right)
\end{equation}
or, from (\ref{eq:08}) we have
\begin{equation}\label{eq:11}
[(I-T)\cd+(T-\Omega)\cl+T\cu]\left(
\begin{array}{c}
x\\y
\end{array}
\right)=\lambda (\cd-\Omega\cl) \left(
\begin{array}{c}
x\\y
\end{array}
\right).
\end{equation}
From (\ref{eq:11}) and (\ref{eq:04}) it follows that
\begin{equation*}
\left(
\begin{array}{ccc}
(1-\tau_1)A&&-\tau_1 B
\\
(\tau_2-\omega_2)B^T&&(1-a\omega_2)Q
\end{array}
\right)\left(
\begin{array}{c}
x\\y
\end{array}
\right)=\lambda \left(
\begin{array}{ccc}
A&&0
\\
-\omega_2B^T&&(1-a\omega_2)Q
\end{array}
\right) \left(
\begin{array}{c}
x\\y
\end{array}
\right).
\end{equation*}
Decoupling we obtain
\begin{equation*}
\begin{cases}
(1-\tau_1)Ax-\tau_1 By=\lambda Ax\\
(\tau_2-\omega_2)B^Tx+(1-a\omega_2)Qy=-\lambda \omega_2 B^T x+\lambda(1-a\omega_2)Qy
\end{cases}
\end{equation*}
or equivalently
\begin{equation}\label{eq:12}
\begin{cases}
(1-\tau_1-\lambda)x=\tau_1 A^{-1}By\\
(1-\lambda)(1-a\omega_2)y=[(1-\lambda) \omega_2-\tau_2] Q^{-1}B^T x.
\end{cases}
\end{equation}
Multiplying the first equality in (\ref{eq:12}) by $Q^{-1}B^T$, we obtain
\begin{equation*}
(1-\tau_1-\lambda)Q^{-1}B^Tx=\tau_1 Q^{-1}B^T A^{-1}By,
\end{equation*}
or, when $\lambda\neq 1-\tau_1$, we have
\begin{equation}\label{eq:13}
Q^{-1}B^Tx=\frac{\tau_1}{1-\tau_1-\lambda} Q^{-1}B^T A^{-1}By.
\end{equation}
From (\ref{eq:13}) and the second equality in (\ref{eq:12}) it follows that
\begin{equation}\label{eq:13'}
(1-\lambda)(1-a\omega_2)(1-\lambda-\tau_1)y=[(1-\lambda)\omega_2-\tau_2]\tau_1 J y.
\end{equation}
\noindent If $\lambda=1-\tau_1\neq0$, we have from (\ref{eq:12}) that $By=0$ and $\tau_1(1-a\omega_2)Qy=(\tau_1\omega_2-\tau_2)B^T x$. It then follows that $y=0$ and $x\in\cn(B^T)$. Hence, $\lambda=1-\tau_1$ is an eigenvalue of $\ch(\tau_1, \tau_2, \omega_2,a)$ with the corresponding eigenvector $(x^T, 0)^T$, where $x\in\cn(B^T)$. Therefore, because of (\ref{eq:13'}), the eigenvalues $\lambda$ (except for $\lambda=1-\tau_1$) of the iteration matrix $\ch(\tau_1,\tau_2,\omega_2,a)$ of the GMESOR($a$) method and the eigenvalues $\mu$ of the matrix J are related through the functional relationship
\begin{equation*}
(1-\lambda)(1-a\omega_2)(1-\lambda-\tau_1)=[(1-\lambda)\omega_2-\tau_2]\tau_1 \mu,
\end{equation*}
namely $\lambda$ satisfies the quadratic equation (\ref{eq:10}). \qed
}
From the above theorem we can obtain the following functional relationships for the GESOR($a$), SOR-like($a$) and GSOR($a$) methods.
\begin{cor} \label{cor: fr112} Under the hypothesis of Theorem \ref{theo:func_rel}\\ \noindent
1. The nonzero eigenvalues of the iteration matrix $\ch(\tau,\omega_2,a)$ of the GESOR($a$) method are given by $\lambda=1-\tau$ or if $a{\omega_2}\neq {1}$ by
\begin{equation}
\lambda^2+\lambda\left(\tau-2+\frac{\tau\omega_2}{1-a\omega_2}\mu\right)+1-\tau+\frac{\tau(\tau-\omega_2)}{1-a\omega_2}\mu=0. \label{eq:15}
\end{equation}
2. The nonzero eigenvalues of the iteration matrix $\cl(\omega,a)$ of the SOR-like($a$) method are given by $\lambda=1-\omega$ or if $a{\omega}\neq {1}$ by
\begin{equation}
\lambda^2+\lambda\left(\omega-2+\frac{\omega^2}{1-a\omega}\mu\right)+1-\omega=0. \label{eq:16}
\end{equation}
3. The nonzero eigenvalues of the iteration matrix $\cl(\omega_1,\omega_2, a)$ of the GSOR($a$) method are given by $\lambda=1-\omega_1$ or if $a{\omega_2}\neq {1}$ by
\begin{equation}
\lambda^2+\lambda\left(\omega_1-2+\frac{\omega_1\omega_2}{1-a\omega_2}\mu\right)+1-\omega_1=0. \label{eq:17}
\end{equation}
\end{cor}
{\bf{Proof}} The iteration matrix $\ch(\tau,\omega_2, a)$ is obtained by letting $\tau=\tau_1=\tau_2$ in $\ch(\tau_1, \tau_2, \omega_2, a)$. By following a similar approach as in the proof of Theorem \ref{theo:func_rel} we find the functional relationship (\ref{eq:15}). Similarly, we find (\ref{eq:16}) and (\ref{eq:17}).\qed 
Note that the above functional relationships are generalizations of the original SOR-like and GSOR methods. Indeed, if $a=0$, then from (\ref{eq:16}) we obtain the functional relationship of the SOR-like method \cite{GolWuYuan2001}, whereas from (\ref{eq:17}) we obtain the functional relationship of the GSOR method \cite{BaiParlWang2005}.\qed 

Another preconditioning matrix $R$, which is formed by the upper triangular part of $\ca$ is the following
\begin{equation}
R=\cd-\Omega \cu.
\label{eq:17b}
\end{equation}
Using (\ref{eq:17b}) in (\ref{eq:06}) then (\ref{eq:05}) becomes the backward form of the GMESOR($a$) method, which will be referred to as the Generalized Modified Extrapolated Backward SOR($a$) (GMEBSOR($a$)) method. From (\ref{eq:06}), because of (\ref{eq:17b}), the iteration matrix of the GMEBSOR($a$) method is given by
\begin{equation*}
\ck{(\tau_1,\tau_2,\omega_1,\omega_2, a)}=I-{(\cd-\Omega\cu)}^{-1}T\ca
\end{equation*}
or
\begin{equation}\label{eq:5_08k}
\ck{(\tau_1,\tau_2,\omega_1,\omega_2, a)}={(\cd-\Omega\cu)}^{-1}[(I-T)\cd+(T-\Omega)\cu+T\cl]
\end{equation}
and
\begin{equation}\label{eq:5_09k}
k{(\tau_1,\tau_2,\omega_1,\omega_2, a)}={(\cd-\Omega\cu)}^{-1}Tb.
\end{equation}
For ${(\cd-\Omega\cu)}^{-1}$ to exist we require
\begin{equation}\label{eq:5_06bbbbu}
\det(\cd-\Omega\cu)\neq 0.
\end{equation}
Because of (\ref{eq:04})
\begin{equation}\label{eq:5_06bbbcu}
R=\cd-\Omega\cu=\left(
\begin{array}{cc}
A & -\omega_1 B\\
0 & [1-(1-a)\omega_2] Q
\end{array}\right).
\end{equation}
Therefore, (\ref{eq:5_06bbbbu}) becomes
\begin{equation*}
\det(\cd-\Omega\cu)=[1-(1-a)\omega_2]^{n}\det{A}\det{Q}\neq 0
\end{equation*}
or
\begin{equation}\label{eq:5_10bu}
(1-a)\omega_2\neq 1
\end{equation}
since the matrix $A$ is symmetric positive definite and the matrix $Q$ is nonsingular.
The GMEBSOR($a$) method has the following algorithmic form.

{\sc The GMEBSOR($a$) Method:}{\em{ Let $Q\in\mathbb{R}^{n\times n}$ be a nonsingular and symmetric matrix.
Given initial vectors $x^{(0)}\in\mathbb{R}^m$ and
$y^{(0)}\in\mathbb{R}^n$, and the parameters $\tau_1,\tau_2\neq 0,  \
\omega_1, \omega_2, a\in\mathbb{R}$ with $(1-a)\omega_2\neq 1$. For $k=0,1,2,...$ until the
iteration sequence
$\{({x^{(k)}}^T, {y^{(k)}}^T)^T\}$ is convergent, compute \\ \newline
\vspace{0.3cm}
$y^{(k+1)}=y^{(k)}+\displaystyle\frac{\tau_2}{1-(1-a)\omega_2}Q^{-1}(B^Tx^{(k)}-b_2)$
\\
\hspace{0.3cm}$x^{(k+1)}=(1-\tau_1)x^{(k)}+ A^{-1}\left\{\tau_1(b_1-By^{(k)})-\omega_1 B(y^{(k+1)}-y^{(k)})\right\} $,\\
\newline
\noindent where Q is an approximate (preconditioning)
matrix of the Schur complement matrix $B^TA^{-1}B$.}}

As a by-product of the GMEBSOR($a$) method we obtain the backward schemes corresponding to the GESOR($a$) and GSOR($a$) methods. For $\tau=\tau_1=\tau_2$, we have the GEBSOR($a$) method and for $\tau_1=\omega_1$ and  $\tau_2=\omega_2$ we have the GBSOR($a$) method.

\begin{cor}
Under the hypothesis of Theorem \ref{theo:func_rel}\\ \noindent
1. The nonzero eigenvalues of the iteration matrix $\ck(\tau_1, \tau_2, \omega_1, \omega_2,a)$ of the GMEBSOR($a$) method are given by $\lambda=1-\tau_1$ or if $(1-a){\omega_2}\neq{1}$ by
\begin{equation}
\lambda^2+\lambda\left(\tau_1-2+\frac{\tau_2\omega_1}{1-(1-a)\omega_2}\mu\right)+1-\tau_1+\frac{\tau_2(\tau_1-\omega_1)}{1-(1-a)\omega_2}\mu=0. \label{eq:18}
\end{equation}
2.  The nonzero eigenvalues of the iteration matrix $\ck(\tau,\omega_1,\omega_2,a)$ of the GEBSOR($a$) method are given by $\lambda=1-\tau$ or if $(1-a){\omega_2}\neq{1}$ by
\begin{equation}
\lambda^2+\lambda\left(\tau-2+\frac{\tau\omega_1}{1-(1-a)\omega_2}\mu\right)+1-\tau+\frac{\tau(\tau-\omega_1)}{1-(1-a)\omega_2}\mu=0. \label{eq:19}
\end{equation}
3.  The nonzero eigenvalues of the iteration matrix $\cm(\omega_1,\omega_2,a)$ of the GBSOR($a$) method are given by $\lambda=1-\omega_1$ or if $(1-a){\omega_2}\neq{1}$ by
\begin{equation}
\lambda^2+\lambda\left(\omega_1-2+\frac{\omega_1\omega_2}{1-(1-a)\omega_2}\mu\right)+1-\omega_1=0. \label{eq:20b}
\end{equation}
\end{cor}
{\bf{Proof}}
Following a similar approach as in the proof of Theorem \ref{theo:func_rel} and using the iteration matrix $\ck(\tau_1, \tau_2,\omega_1, \omega_2, a)$ given by (\ref{eq:5_08k}), we find the functional relationship (\ref{eq:18}). Similarly, we find the functional relationships (\ref{eq:19}) and (\ref{eq:20b}).\qed

Note that the GMEBSOR($a$) method has four parameters instead of three as the GMESOR($a$) method whereas the GEBSOR($a$) method has three parameters instead of two as the GESOR($a$) method. If $a=1$, then (\ref{eq:20b}) becomes the functional relationship of the GSOR($a$) method.
}

\subsection{Convergence}
\label{sec:3}
In this section we develop the convergence analysis of the { {GSOR($a$) and}} GMESOR method{ {s as well as their corresponding backward counterparts.}}
In particular, we derive sufficient conditions for the { {GSOR($a$) and the}} GMESOR method to converge under the assumption that the eigenvalues of the matrix $J$ are all real. The sign of $J$'s eigenvalues depends upon the properties of the matrix $Q$. Specifically, we assume that the matrix $Q$ is symmetric positive or negative definite.
{ {
\subsubsection{The GSOR($a$) method}
The next theorem provides sufficient conditions for the GSOR($a$) method to converge if the matrix $Q$ is symmetric positive definite and $a\neq 0$.
\begin{theo}\label{theo:conv_gsor}
Let $A\in\mathbb{R}^{m\times m}$ and $Q\in\mathbb{R}^{n\times n}$ be
symmetric positive definite and $B\in\mathbb{R}^{m\times n}$ be of
full column rank. Denote the minimum and the maximum eigenvalues of the matrix $J=Q^{-1}B^TA^{-1}B$ by $\mu_{min}$ and $\mu_{max}$, respectively. Then $\rho({\cl}{(\omega_1,\omega_2, a)})<1$ if
the parameters $\omega_1\; \mbox{and}\; \omega_2$ lie in any
case of Table \ref{tab1}.
\begin{table}[htbp]
{ {\caption{Sufficient conditions for the
 GSOR($a$) method to converge if $\mmin>0$. } \label{tab1}
\begin{tabular}{|c|c|c|c|c|c|c|c|c|c|c|}
\hline
   ${\mbox{Condition}}$&$\mbox{Cases}$ & $\omega_2 - \mbox{Domain}$
   & $\omega_1 - \mbox{Domain}$\\
    \hline\hline
 $a>0$&1
 & $0<\omega_2<\displaystyle\frac{2(2-\omega_1)}{\omega_1\mmax+2a(2-\omega_1)}$&$0<\omega_1<2 $
 \\ \hline &2& $\omega_2<\displaystyle\frac{2(2-\omega_1)}{\omega_1\mmax+2a(2-\omega_1)}(<0)$ & $0<\omega_1<\displaystyle\frac{4a}{2a-\mmax}$
  \\ \cline{2-2} \cline{3-4} $a<0$&3
 &   $0<\omega_2$& $0<\omega_1<\displaystyle\frac{4a}{2a-\mmax}$
  \\ \cline{2-2} \cline{3-4} &4
 & $0<\omega_2<\displaystyle\frac{2(2-\omega_1)}{\omega_1\mmax+2a(2-\omega_1)}$ &$\displaystyle\frac{4a}{2a-\mmin}<\omega_1<2$
 \\ \hline
\end{tabular}}}
\end{table}
\end{theo}
{\bf{Proof}}  Recall (Corollary \ref{cor: fr112}) that $\lambda=1-\omega_1 \neq 0$ is an eigenvalue of $\cl{(\omega_1,\omega_2, a)}$ and if $\lambda\neq 1-\omega_1$ then the eigenvalues of $\cl{(\omega_1,\omega_2, a)}$ are given by (\ref{eq:17}). If $\lambda=1-\omega_1 \neq 0$, then the GSOR($a$) method is convergent if and only if $|\lambda|<1$, that is $|1-\omega_1|<1$, or
\begin{equation}\label{eq:31gsor}
0<\omega_1<2.
\end{equation}
If $\lambda\neq 1-\omega_1$ and $a{\omega_2}\neq{1}$, then
(\ref{eq:17}) holds and by Lemma 2.1 page 171 of \cite{Young} it follows that the GSOR($a$) method is convergent if and only if
\begin{equation}\label{eq:3g2}
|c|<1 \;\; \mbox{and} \;\; |b|<1+c
\end{equation}
where
\begin{equation}\label{eq:3g3}
c=1-\omega_1
\end{equation}
and
\begin{equation}\label{eq:3g4}
b=2-\omega_1-\frac{\omega_1\omega_2\mu}{1-a\omega_2}.
\end{equation}
From the first part of (\ref{eq:3g2}), because of (\ref{eq:3g3}), it follows that (\ref{eq:31gsor}) holds also in this case. From the second part of (\ref{eq:3g2}), because of (\ref{eq:3g3}) and (\ref{eq:3g4}), it follows that
\begin{equation*}
\left|2-\omega_1-\frac{\omega_1\omega_2\mu}{1-a\omega_2}\right|<2-\omega_1
\end{equation*}
or
\begin{equation}\label{eq:3g5}
0<\frac{\omega_2}{1-a\omega_2}<\frac{2(2-\omega_1)}{\omega_1\mu}.
\end{equation}
In order for (\ref{eq:3g5}) to hold we distinguish two cases. Case I: $\omega_2>0$ and $1-a\omega_2>0$ and Case II: $\omega_2<0$ and $1-a\omega_2<0$. For each of theses cases we will distinguish two subcases. (i) $a>0$ (ii) $a<0$.
In the sequel we will study the subcase (i) of Case I, since the other cases can be treated similarly.
For subcase (i) of Case I
\begin{equation}\label{eq:3g6}
0<\omega_2<\frac{1}{a}.
\end{equation}
From (\ref{eq:3g5}), we have
\begin{equation}\label{eq:3g7}
[\omega_1\mu+2a(2-\omega_1)]\omega_2<2(2-\omega_1).
\end{equation}
We distinguish two subcases: $(\mbox{i}_1)\;\; \omega_1\mu+2a(2-\omega_1)>0$ and $(\mbox{i}_2)\;\; \omega_1\mu+2a(2-\omega_1)<0$. In the sequel we will only treat subcase $(\mbox{i}_1)$ since the other case can be treated similarly. If $\omega_1\mu+2a(2-\omega_1)>0$ then
\begin{equation}\label{eq:3g8}
4a>\omega_1(2a-\mu).
\end{equation}
Next, we distinguish three subcases: (a) $a\geq \frac{1}{2}\mmax$ (b) $a\leq \frac{1}{2}\mmin$ (c) $\frac{1}{2}\mmin<a<\frac{1}{2}\mmax$.\\ \noindent
(a) $a\geq \frac{1}{2}\mmax$. From (\ref{eq:3g8}) we have
\begin{equation}\label{eq:3gs1}
\omega_1<\frac{4a}{2a-\mmin}.
\end{equation}
Combining (\ref{eq:31gsor}) and (\ref{eq:3gs1}), it follows that
\begin{equation*}
0<\omega_1<min\left\{2,\frac{4a}{2a-\mmin}\right\},
\end{equation*}
or
\begin{equation}\label{eq:3gs2}
0<\omega_1<2.
\end{equation}
Moreover, from (\ref{eq:3g7}) we have
\begin{equation}\label{eq:3gs3}
0<\omega_2<\frac{2(2-\omega_1)}{\omega_1\mmax+2a(2-\omega_1)}
\end{equation}
which, because of (\ref{eq:3g6}), becomes
\begin{equation}\label{eq:3gs4}
0<\omega_2<min\left\{\frac{1}{a},\frac{2(2-\omega_1)}{\omega_1\mmax+2a(2-\omega_1)}\right\},
\end{equation}
which yields (\ref{eq:3gs3}) again. Therefore, for case (a) we have that (\ref{eq:3gs2}) and (\ref{eq:3gs3}) hold.\\ \noindent
(b) $a\leq \frac{1}{2}\mmin$. From (\ref{eq:3g8}) we have
\begin{equation}\label{eq:3gs12}
\omega_1>\frac{4a}{2a-\mmax}.
\end{equation}
Combining (\ref{eq:31gsor}) and (\ref{eq:3gs12}), it follows that
\begin{equation*}
max\left\{0,\frac{4a}{2a-\mmax}\right\}<\omega_1<2,
\end{equation*}
which yields (\ref{eq:3gs2}).
Therefore, for case (b) we have that (\ref{eq:3gs2}) and (\ref{eq:3gs3}) hold also as in case (a).\\ \noindent
(c) $\frac{1}{2}\mmin<a<\frac{1}{2}\mmax$.
Let $\alpha, \beta$  be two positive integers such that $\mu_{\alpha}=\max\{\mu | \mu\leq 2\alpha\}$, $\mu_{\beta}=\min\{\mu | \mu\geq 2\alpha\}$.
Next, we distinguish two cases: {(i)} $\mu_{min}\leq \mu\leq \mu_{\alpha}$, {(ii)} $\mu_{\beta} \leq \mu\leq \mu_{max}$.\\
Case {(i)}: $\mu_{min}\leq \mu\leq \mu_{\alpha}$. Following a similar approach as in Case ({a}), we have that
(\ref{eq:3gs2}) holds and
\begin{equation}\label{eq:5_3gs3aa}
0<\omega_2<\frac{2(2-\omega_1)}{\omega_1\mu_\alpha+2a(2-\omega_1)}.
\end{equation}
Case {(ii)}: $\mu_{\beta} \leq \mu\leq \mu_{max}$. Following a similar approach as in Case ({b}), we have that
(\ref{eq:3gs2}) holds and
\begin{equation}\label{eq:5_3gs3bb}
0<\omega_2<\frac{2(2-\omega_1)}{\omega_1\mmax+2a(2-\omega_1)}.
\end{equation}
Combining (\ref{eq:5_3gs3aa}) and (\ref{eq:5_3gs3bb}) it follows that
\begin{equation}\label{eq:5_3gs3cc}
0<\omega_2<\min\left\{\frac{2(2-\omega_1)}{\omega_1\mu_\alpha+2a(2-\omega_1)}, \frac{2(2-\omega_1)}{\omega_1\mmax+2a(2-\omega_1)}\right\},
\end{equation}
which is equivalent to (\ref{eq:3gs3}).
Hence, case 1 of table \ref{tab1} is proved.
Following a similar treatment we can prove the rest of the cases of Table \ref{tab1}.
\qed

\begin{cor}\label{cor:conv03b}
Under the hypothesis of Theorem \ref{theo:conv_gsor} and if $a=0$ then ${\rho(\cl{(\omega_1, \omega_2)})<1}$ if
\begin{equation}
    0<\omega_1<2 \ \ \mbox{and}\ \
   0<\omega_2<\frac{2(2-\omega_1)}{\omega_1\mu_{max}}.\label{eq:43}
\end{equation}
\end{cor}
{\bf{Proof}}
If we let $a=0$ in (\ref{eq:17}) and follow a similar approach as in the proof of Theorem
\ref{theo:conv_gsor} we can verify that (\ref{eq:43}) holds. \qed
Note that (\ref{eq:43}) was also obtained in \cite{BaiParlWang2005}.
The following corollary gives sufficient conditions for the GBSOR($a$) method to converge.
\begin{cor}\label{cor:conv3}
Under the hypothesis of Theorem \ref{theo:conv_gsor}, $\rho(\cm{(\omega_1, \omega_2,a)})<1$ if the parameters $\omega_1$ and $\omega_2$ lie in any case of Table \ref{tab2}.
\begin{table}[htbp]
{ {\caption{Sufficient conditions for the
 GBSOR($a$) method to converge  if $\mmin>0$. } \label{tab2}
\begin{tabular}{|c|c|c|c|c|c|c|c|c|c|c|}
\hline
   ${\mbox{Condition}}$&$\mbox{Cases}$ & $\omega_2 - \mbox{Domain}$
   & $\omega_1 - \mbox{Domain}$\\
    \hline\hline
 $a<1$&1
 & $0<\omega_2<\displaystyle\frac{2(2-\omega_1)}{\omega_1\mmax+2(1-a)(2-\omega_1)}$&$0<\omega_1<2 $
 \\ \hline &2& $\omega_2<\displaystyle\frac{2(2-\omega_1)}{\omega_1\mmax+2(1-a)(2-\omega_1)}<0$ & $0<\omega_1<\displaystyle\frac{4(1-a)}{2(1-a)-\mmax}$
  \\ \cline{2-2} \cline{3-4} $1<a$&3
 &   $0<\omega_2$& $0<\omega_1<\displaystyle\frac{4(1-a)}{2(1-a)-\mmax}$
  \\ \cline{2-2} \cline{3-4} &4
 & $0<\omega_2<\displaystyle\frac{2(2-\omega_1)}{\omega_1\mmax+2(1-a)(2-\omega_1)}$ &$\displaystyle\frac{4(1-a)}{2(1-a)-\mmin}<\omega_1<2$
 \\ \hline
\end{tabular}}}
\end{table}
\end{cor}
{\bf{Proof}} Using the functional relationship (\ref{eq:20b}) and
following a similar approach as in the proof of Theorem \ref{theo:conv_gsor} we have
\begin{equation}\label{eq:cor8_12}
0<\omega_1<2\; \mbox{and}\; 0<\frac{\omega_2}{1-(1-a)\omega_2}<\frac{2(2-\omega_1)}{\omega_1\mu}.
\end{equation}
Note that the second part of (\ref{eq:cor8_12}) is the same as (\ref{eq:3g5}) where now $1-a$ appears instead of $a$. This occurs because the preconditioning matrix $R$ is given by (\ref{eq:17b}) and $\mathcal{U}$ is expressed in (\ref{eq:04}) in terms of $1-a$. Therefore, if we let $1-a$ in place of $a$ in Table \ref{tab1}, we obtain Table \ref{tab2}. \qed
\begin{cor}\label{cor:conv333}
Under the hypothesis of corollary \ref{cor:conv3} and if $a=1$ then ${\rho(\cm{(\omega_1, \omega_2,1)})<1}$ if
\begin{equation}
    0<\omega_1<2 \ \ \mbox{and}\ \
   0<\omega_2<\frac{2(2-\omega_1)}{\omega_1\mu_{max}}.\label{eq:4333}
\end{equation}
\end{cor}
{\bf{Proof}} If we let $a=1$ in  (\ref{eq:cor8_12}) then (\ref{eq:4333}) follows immediately. \qed
If the matrix Q is symmetric negative definite, then we have the following theorem.
\begin{theo}\label{theo:conv_gsornon}
Let $A\in\mathbb{R}^{m\times m}$  be
symmetric positive definite, $B\in\mathbb{R}^{m\times n}$ be of
full column rank and $Q\in\mathbb{R}^{n\times n}$ be
symmetric negative definite. Denote the minimum and the maximum
eigenvalues of the matrix $J=Q^{-1}B^TA^{-1}B$ by $\mu_{min}$ and $\mu_{max}$, respectively.
Then $\rho({\cl}{(\omega_1,\omega_2,a)})<1$ if the
parameters $\omega_1\; \mbox{and}\; \omega_2$ lie in the following cases of Table \ref{tab1b}.
\begin{table}[htbp]
{ {\caption{Sufficient conditions for the
GSOR($a$) method to converge if $\mmax<0$.} \label{tab1b}
\begin{tabular}{|c|c|c|c|c|c|c|c|c|c|c|}
\hline
   ${\mbox{Condition}}$&$\mbox{Cases}$ & $\omega_2 - \mbox{Domain}$
   & $\omega_1 - \mbox{Domain}$\\
    \hline\hline
 $a<0$&1
 & $\displaystyle\frac{2(2-\omega_1)}{\omega_1\mmin+2a(2-\omega_1)}<\omega_2<0$&$0<\omega_1<2 $
 \\ \hline &2& $(0<)\displaystyle\frac{2(2-\omega_1)}{\omega_1\mmin+2a(2-\omega_1)}<\omega_2$ & $0<\omega_1<\displaystyle\frac{4a}{2a-\mmin}$
  \\ \cline{2-2} \cline{3-4} $a>0$&3
 &   $\omega_2<0$& $0<\omega_1<\displaystyle\frac{4a}{2a-\mmin}$
  \\ \cline{2-2} \cline{3-4} &4
 & $\displaystyle\frac{2(2-\omega_1)}{\omega_1\mmin+2a(2-\omega_1)}<\omega_2<0$ &$\displaystyle\frac{4a}{2a-\mmax}<\omega_1<2$
 \\ \hline
\end{tabular}}}
\end{table}
\end{theo}
{\bf{Proof}} Using the functional relationship (\ref{eq:17}) and
following a similar approach as in the proof of Theorem \ref{theo:conv_gsor} taking into consideration that ${\mmax<0}$ we have
\begin{equation}\label{eq:th32}
0<\omega_1<2\;\;\;\mbox{and}\;\;\; \frac{2(2-\omega_1)}{\omega_1\mu}<\frac{\omega_2}{1-a\omega_2}<0.
\end{equation}
From (\ref{eq:th32}) the cases presented in Table \ref{tab1b} can be readily verified.\qed

\begin{cor}\label{cor:conv4}
Under the hypothesis of Theorem \ref{theo:conv_gsornon} and  if $a=0$ then ${\rho(\cl{(\omega_1,\omega_2)})<1}$ if
\begin{equation}
    0<\omega_1<2\ \ \mbox{and}\ \
    \frac{2(2-\omega_1)}{\omega_1\mu_{min}}<\omega_2<0.\label{eq:44}
\end{equation}
\end{cor}
{\bf{Proof}} Using the functional relationship (\ref{eq:17}) and
following the proof of Theorem \ref{theo:conv_gsornon} we have that if $a=0$ in (\ref{eq:th32}) then (\ref{eq:44}) follows. \qed
The above result was also obtained in \cite{BaiParlWang2005}.
\begin{theo}\label{theo:conv_gbsornon}
Let $A\in\mathbb{R}^{m\times m}$  be
symmetric positive definite, $B\in\mathbb{R}^{m\times n}$ be of
full column rank and $Q\in\mathbb{R}^{n\times n}$ be
symmetric negative definite. Denote the minimum and the maximum
eigenvalues of the matrix $J=Q^{-1}B^TA^{-1}B$ by $\mu_{min}$ and $\mu_{max}$, respectively.
Then $\rho({\cm}{(\omega_1,\omega_2,a)})<1$ if the
parameters $\omega_1\; \mbox{and}\; \omega_2$ lie in the following cases of Table \ref{tab1b}.
\begin{table}[htbp]
{ {\caption{Sufficient conditions for the
GBSOR($a$) method to converge if $\mmax<0$.} \label{tab1b4}
\begin{tabular}{|c|c|c|c|c|c|c|c|c|c|c|}
\hline
   ${\mbox{Condition}}$&$\mbox{Cases}$ & $\omega_2 - \mbox{Domain}$
   & $\omega_1 - \mbox{Domain}$\\
    \hline\hline
 $a<1$&1
 & $\displaystyle\frac{2(2-\omega_1)}{\omega_1\mmin+2(1-a)(2-\omega_1)}<\omega_2<0$&$0<\omega_1<2 $
 \\ \hline &2& $(0<)\displaystyle\frac{2(2-\omega_1)}{\omega_1\mmin+2(1-a)(2-\omega_1)}<\omega_2$ & $0<\omega_1<\displaystyle\frac{4(1-a)}{2(1-a)-\mmin}$
  \\ \cline{2-2} \cline{3-4} $a>1$&3
 &   $\omega_2<0$& $0<\omega_1<\displaystyle\frac{4(1-a)}{2(1-a)-\mmin}$
  \\ \cline{2-2} \cline{3-4} &4
 & $\displaystyle\frac{2(2-\omega_1)}{\omega_1\mmin+2(1-a)(2-\omega_1)}<\omega_2<0$ &$\displaystyle\frac{4(1-a)}{2(1-a)-\mmax}<\omega_1<2$
 \\ \hline
\end{tabular}}}
\end{table}
\end{theo}
{\bf{Proof}} Using the functional relationship (\ref{eq:20b}) and
following a similar approach as in the proof of Theorem \ref{theo:conv_gsor} taking into consideration that ${\mmax<0}$ we have
\begin{equation}\label{eq:th324}
0<\omega_1<2\;\;\;\mbox{and}\;\;\; \frac{2(2-\omega_1)}{\omega_1\mu}<\frac{\omega_2}{1-(1-a)\omega_2}<0.
\end{equation}
From (\ref{eq:th324}) the cases presented in Table \ref{tab1b4} can be readily verified.\qed
\begin{cor}\label{cor:conv4b}
Under the hypothesis of Theorem \ref{theo:conv_gbsornon} and  if $a=1$ then ${\rho(\cm{(\omega_1,\omega_2)})<1}$ if
\begin{equation}
    0<\omega_1<2\ \ \mbox{and}\ \
    \frac{2(2-\omega_1)}{\omega_1\mu_{min}}<\omega_2<0.\label{eq:444}
\end{equation}
\end{cor}
{\bf{Proof}} Using the functional relationship (\ref{eq:20b}) and
following the proof of Theorem \ref{theo:conv_gbsornon} we have that if $a=1$ in (\ref{eq:th324}) then (\ref{eq:444}) follows. \qed
In the sequel we study the convergence analysis of the GMESOR method under the same assumptions.

}}

\subsubsection{The GMESOR method}
The next theorem provides sufficient conditions for the GMESOR method to converge if the matrix Q is symmetric positive definite and $a=0$. The study of the case $a\neq 0$ follows a similar but cumbersome approach as it requires many cases to be examined. This study will not have any substantial contribution since the minimum value of the spectral radius of the GMESOR($a$) method is independent of $a$ (Theorem \ref{theo:opt_1}), meaning that for, say $a=0$, the GMESOR method will attain the maximum rate of convergence. So, we are interested to find the convergence ranges of the parameters of the GMESOR($a$) method for the simplified case when $a=0$.
\begin{theo}\label{theo:conv1}
Consider the GMESOR method. Let $A\in\mathbb{R}^{m\times m}$ and $Q\in\mathbb{R}^{n\times n}$ be
symmetric positive definite and $B\in\mathbb{R}^{m\times n}$ be of
full column rank. Denote the minimum and the maximum
eigenvalues of the matrix $J=Q^{-1}B^TA^{-1}B$ by $\mu_{min}$ and $\mu_{max}$, respectively.
Then $\rho({\ch}{(\tau_1,\tau_2,\omega_2)})<1$ if
{{
\begin{equation}
    0<\tau_1<2,\ \ 0<\tau_2<\bar\tau_2{(\mu_{max})} \ \ \mbox{and}\ \
   \underline\omega_2({\mu_{max}})<\omega_2< \bar\omega_2({\mu_{max}})
    \label{eq:30}
\end{equation}
where
\begin{equation}
    \bar\tau_2{(\mu_{max})}=\frac{4}{\tau_1\mu_{max}}, \ \
   \underline\omega_2({\mu_{max}})= \tau_2-\frac{1}{\mu_{max}}\ \ \mbox{and}\ \ \bar\omega_2({\mu_{max}})=\frac{2-\tau_1}{\tau_1\mu_{max}}+\frac{\tau_2}{2}.
    \label{eq:300}
\end{equation}
}}
\end{theo}
{ {{\bf{Proof}} Recall that $\lambda=1-\tau_1 \neq 0$ is an eigenvalue of $\ch{(\tau_1,\tau_2,\omega_2)}$ and if $\lambda\neq 1-\tau_1$ then the eigenvalues of $\ch{(\tau_1,\tau_2,\omega_2)}$ are given by (\ref{eq:10}) where $a=0$. If $\lambda=1-\tau_1 \neq 0$, then the GMESOR method is convergent if and only if $|\lambda|<1$, that is $|1-\tau_1|<1$, or
\begin{equation}\label{eq:31}
0<\tau_1<2,
\end{equation}
which is the first inequality of (\ref{eq:30}).
If $\lambda\neq 1-\tau_1$, then (\ref{eq:10}) holds and by Lemma 2.1 page 171 of \cite{Young}, it follows that the GMESOR method is convergent if and only if (\ref{eq:3g2}) holds
where
\begin{equation}\label{eq:33}
c=1-\tau_1+\tau_1(\tau_2-\omega_2)\mu
\end{equation}
and
\begin{equation}\label{eq:34}
b=2-\tau_1-\tau_1\omega_2\mu.
\end{equation}
From the first inequality of (\ref{eq:3g2}) it follows that
\begin{equation}\label{eq:35}
0<1+c<2.
\end{equation}
From the second inequality of (\ref{eq:3g2}), because of (\ref{eq:33}) and (\ref{eq:34}), we have
\begin{equation*}
|1+c-\tau_1\tau_2\mu|<1+c
\end{equation*}
or
\begin{equation}\label{eq:36}
0<\frac{\tau_1\tau_2\mu}{2}<1+c.
\end{equation}
Combining (\ref{eq:35}) and (\ref{eq:36}), it follows that
\begin{equation}\label{eq:37}
0<\frac{\tau_1\tau_2\mu}{2}<1+c<2.
\end{equation}
In order for (\ref{eq:37}) to hold we must have
\begin{equation*}
0<\frac{\tau_1\tau_2\mu}{2}<2
\end{equation*}
or, because of (\ref{eq:31}),
\begin{equation}\label{eq:38}
0<\tau_2<\frac{4}{\tau_1\mu},
\end{equation}
which proves the second inequality of (\ref{eq:30}).
Inequality (\ref{eq:37}), because of (\ref{eq:33}), becomes
\begin{equation*}
\frac{\tau_1\tau_2\mu}{2}<2-\tau_1+\tau_1\tau_2\mu-\tau_1\omega_2\mu<2
\end{equation*}
which is equivalent to
\begin{equation}\label{eq:39}
\tau_2-\frac{1}{\mu}<\omega_2<\frac{2-\tau_1}{\tau_1\mu}+\frac{\tau_2}{2}.
\end{equation}
By studying  the monotonicity of the right and left hand side of (\ref{eq:39}) with respect to $\mu$ we obtain the third inequality of (\ref{eq:30}).\qed
The convergence conditions for GESOR are given by the following corollary.
{ {\begin{cor}\label{cor:conv1}
Consider the GESOR method. Let $A\in\mathbb{R}^{m\times m}$ and $Q\in\mathbb{R}^{n\times n}$ be
symmetric positive definite and $B\in\mathbb{R}^{m\times n}$ be of
full column rank. Denote the minimum and the maximum
eigenvalues of the matrix $J=Q^{-1}B^TA^{-1}B$ by $\mu_{min}$ and $\mu_{max}$, respectively. Then $\rho({\ch}{(\tau,\omega_2)})<1$ if
\begin{equation}\label{eq:41}
    0<\tau<\bar\tau{(\mu_{max})} \ \ \mbox{and}\ \
    \underline\omega_2(\tau)<\omega_2<\bar\omega_2(\tau),
\end{equation}
where
\begin{equation}\label{eq:0041}
\bar\tau{(\mu_{max})}=
\begin{cases}
\hspace{0.8cm} 2, & \mmax\leq 1\\
\frac{2}{\sqrt{\mmax}}, & \mmax>1,
\end{cases}
\end{equation}
\begin{equation}\label{eq:041}
       \underline\omega_2(\tau)= \tau-\frac{1}{\mu_{max}}\ \ \mbox{and}\ \ \bar\omega_2(\tau)=\frac{2-\tau}{\tau\mu_{max}}+\frac{\tau}{2}.     \end{equation}
\end{cor}
{\bf{Proof}}
Letting $\tau=\tau_1=\tau_2$ in (\ref{eq:30}) we obtain (\ref{eq:41}). \qed
The convergence area for the GESOR method is illustrated in figure 1. Note that as $\mmax$ increases the point of intersection of the two curves $\bar\omega_2(\tau)$ and $\underline\omega_2(\tau)$ moves towards zero and the convergence area of the GESOR method shrinks. However, in practice $\mmax$ usually is $<1$.
\begin{figure}[ht]
\begin{center}
   \scalebox{0.8}{\includegraphics{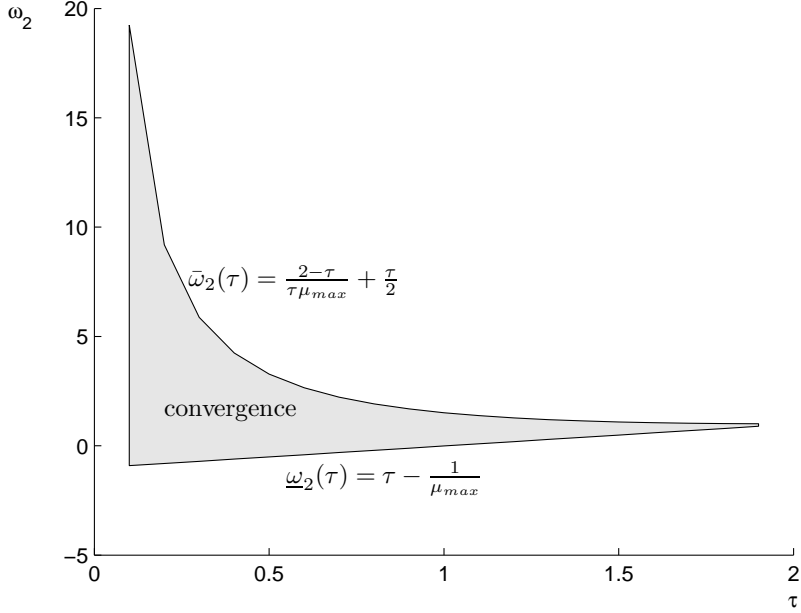}}
  \end{center}\label{fig:5-099}
{ { \caption{Convergence area of the GESOR method for $\mmax=0.99$.}}}
\end{figure}
}}

If the matrix Q is symmetric negative definite and $a=0$ then we have the following theorem.

{ {\begin{theo}\label{theo:conv00}
Consider the GMESOR method. Let $A\in\mathbb{R}^{m\times m}$  be
symmetric positive definite, $B\in\mathbb{R}^{m\times n}$ be of
full column rank and $Q\in\mathbb{R}^{n\times n}$ be
symmetric negative definite. Denote the minimum and the maximum
eigenvalues of the matrix $J=Q^{-1}B^TA^{-1}B$ by $\mu_{min}$ and $\mu_{max}$, respectively. Then $\rho({\ch}{(\tau_1, \tau_2,\omega_2)})<1$ if
{ {
\begin{equation}
    0<\tau_1<2,\ \ \underline\tau_2{(\mu_{min})}<\tau_2<0\ \ \mbox{and}\ \
    \underline\omega_2(\mmin)<\omega_2<\bar\omega_2(\mmin).
    \label{eq:40}
\end{equation}
where
\begin{equation}
   \underline\tau_2{(\mu_{min})}=\frac{4}{\tau_1\mu_{min}},\ \
    \underline\omega_2(\mmin)=\frac{2-\tau_1}{\tau_1\mu_{min}}+\frac{\tau_2}{2} \ \ \mbox{and}\ \   \bar\omega_2(\mmin)=\tau_2-\frac{1}{\mu_{min}}.
    \label{eq:400}
\end{equation}
}}
\end{theo}
{\bf{Proof}} { {Following a similar approach as in the proof of Theorem \ref{theo:conv1} and using the functional relationship (\ref{eq:10}) we can prove (\ref{eq:40}).}}\qed
%
%
\begin{cor}\label{cor:conv2}
Consider the GESOR method. Let $A\in\mathbb{R}^{m\times m}$  be
symmetric positive definite, $B\in\mathbb{R}^{m\times n}$ be of
full column rank and $Q\in\mathbb{R}^{n\times n}$ be
symmetric negative definite. Denote the minimum and the maximum
eigenvalues of the matrix $J=Q^{-1}B^TA^{-1}B$ by $\mu_{min}$ and $\mu_{max}$, respectively. Then $\rho({\ch}{(\tau, \omega_2)})<1$ if
{ {\begin{equation}
    \underline\tau_1{(\mu_{min})}<\tau<2,\ \
    \underline\omega_2(\mmin)<\omega_2<\bar\omega_2(\mmin)\ \ \mbox{and}\ \ \mmin>1,
    \label{eq:42}
\end{equation}
where
\begin{equation}
    \underline\tau_1{(\mu_{min})}=\frac{2}{\sqrt{\mmin}},\ \
    \underline\omega_2(\mmin)=\frac{2-\tau}{\tau\mu_{min}}+\frac{\tau}{2}\ \ \mbox{and}\ \        \bar\omega_2(\mmin)=\tau-\frac{1}{\mu_{min}}.\label{eq:420}
\end{equation}
}}
\end{cor}
{\bf{Proof}}
(\ref{eq:42}) is proved by following a similar approach as in the proof of Theorem \ref{theo:conv1} and using the functional relationship (\ref{eq:15}).\qed

}}
}}

\noindent

 \noindent

\subsection{Optimum parameters}
In this section we determine optimum values for the parameters of the iterative methods studied in the present section under the hypothesis that $a\neq 0$ and the eigenvalues of the matrix $J$ are real. We assume that $Q$ is a symmetric positive or negative definite matrix.

{ {\subsubsection{The GSOR($a$) method}
In the following theorem the optimum parameters for the GSOR($a$) method are determined assuming that the matrix $Q$ is symmetric positive definite.
\begin{theo}\label{theo:opt_gsor}
Consider the GSOR($a$) method. Let $A\in\mathbb{R}^{m\times m}$ and $Q\in\mathbb{R}^{n\times n}$ be symmetric positive definite and $B\in\mathbb{R}^{m\times n}$ be of
full column rank. Denote the minimum and the maximum eigenvalues of the matrix $J=Q^{-1}B^TA^{-1}B$ by
$\mu_{min}$ and $\mu_{max}$, respectively.
Then the spectral radius of the GSOR($a$) method, $\rho(\cl({\omega_1, \omega_2,a}))$, is minimized for any $a\neq {-\sqrt{\mmin\mmax}}$ at
\begin{equation}\label{eq:3u2ux}
\omega_{1_{opt}}=\frac{4\sqrt{\mu_{min}\mu_{max}}}{(\sqrt{\mu_{min}}+\sqrt{\mu_{max}})^2}\ \ \mbox{and} \ \ \omega_{2_{opt}}=\frac{1}{a+\sqrt{\mu_{min}\mu_{max}}}
\end{equation}
and its corresponding value is
\begin{equation}\label{eq:4u2ux}
\rho(\cl{(\omega_{1_{opt}},\omega_{2_{opt}},a)})=\frac{\sqrt{\mu_{max}}-\sqrt{\mu_{min}}}{\sqrt{\mu_{max}}+\sqrt{\mu_{min}}}.
\end{equation}
\end{theo}
{\bf{Proof}} The functional relationship (\ref{eq:17})  may be written as follows
\begin{equation}\label{eq:5_s01bbxx}
(\lambda+\omega_1-1)(\lambda-1)=-\lambda\omega_1\hat\omega_2\mu
\end{equation}
\noindent where
\begin{equation}\label{eq:5_s01bbxxaxax}
\hat\omega_2=\frac{\omega_2}{1-a\omega_2}
\end{equation}
with $a\omega_2\neq 1$. The optimum values of $\omega_1\;\;\mbox{and}\;\; \hat\omega_2$
will be determined such that
\begin{equation}\label{eq:5_s02bbxx}
\rho{(\cl(\omega_1,\omega_2, a))}=\max_{\mu_{min}\leq\mu\leq\mu_{max}}|\lambda|
\end{equation}
is minimum. The real roots of (\ref{eq:5_s01bbxx}) are the intersection points
of the parabola
\begin{equation}\label{eq:5_s03bbxx}
g_{\omega_1}(\lambda)=\frac{(\lambda+\omega_1-1)(\lambda-1)}{\omega_1\hat\omega_2}
\end{equation}
and the straight lines
\begin{equation}\label{eq:5_s04bbxx}
h(\lambda)=-\lambda\mu, \;\; 0<\mu_{min}\leq\mu\leq\mu_{max}.
\end{equation}
Following a similar argument as in \cite{Var1} page 111,
$h(\lambda)$ are straight lines through the point $(0,0)$
and $g_{\omega_1}(\lambda)$ is a parabola passing through the point (1,0).
The discriminant of (\ref{eq:17}) is
\begin{equation}\label{eq:5_s05bbxxa}
\Delta(\omega_1,\hat\omega_2,\mu)=(2-\omega_1-\omega_1\hat\omega_2\mu)^2-4(1-\omega_1).
\end{equation}
Note that $\Delta(\omega_1,\hat\omega_2,\mu)\leq 0$ for $0<\omega_1\leq \tilde\omega_1(\mu)$ and
$\Delta(\omega_1,\hat\omega_2,\mu)\geq 0$ for ${\tilde\omega_1(\mu)\leq\omega_1<2}$, where
\begin{equation}\label{eq:5_mar}
\tilde\omega_1(\mu)=\frac{4\hat\omega_2\mu}{(1+\hat\omega_2\mu)^2}.
\end{equation}
If $0<\omega_1\leq \tilde\omega_1(\mu)$  then the value of $\rho(\cl(\omega_1, \omega_2,a))$ is
\begin{equation}\label{eq:5_s2bbbxxa}
|\tilde\lambda_1|=|\tilde\lambda_N|=(1-\omega_1)^{1/2},
\end{equation}
where $\tilde\lambda_1$ and $\tilde\lambda_N$ are the two conjugate complex roots of (\ref{eq:17}).
Furthermore, (\ref{eq:5_s2bbbxxa}) is a decreasing function of $\omega_1$.
In case $\tilde\omega_1(\mu)\leq\omega_1<2$ the roots of (\ref{eq:17}) can be geometrically interpreted as the intersection of the curves $g_{\omega_1}(\lambda)$ and $h(\lambda)=-\lambda\mu$, as illustrated in figure 2, where we have assumed, without loss of generality, that $h(\lambda) \equiv h_1(\lambda)=-\lambda\mu_{max}$.
\begin{figure}[ht]
  \begin{center}
   \scalebox{0.5}{\includegraphics{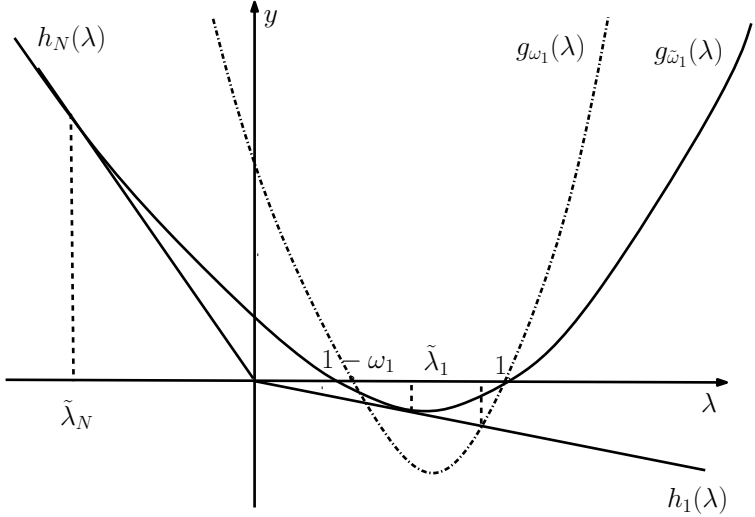}}
  \end{center}\label{fig:5-01}
  { {\caption{Graphs of $g_{\omega_1}(\lambda), h_1(\lambda)$ and $h_N(\lambda)$ in case the roots of (\ref{eq:5_s01bbxx}) are real.}}}
\end{figure}
\noindent The largest abscissa of the two points of
intersection of $h(\lambda)$ and $g_{\omega_1}(\lambda)$ decreases with increasing $\omega_1$. Indeed as $\omega_1$ increases, the intersection point $(1-\omega_1,0)$ of
$g_{\omega_1}(\lambda)$ with the $O\lambda$ axis is moving towards to zero
until $g_{\omega_1}(\lambda)$ becomes tangent to $h(\lambda)$. Thus, for the fixed eigenvalue $\mu$ of $J$, the value of $\omega_1$ which minimizes the zero of largest modulus of (\ref{eq:17}) is $\tilde\omega_1(\mu)$.
Note that the straight lines
$h_1(\lambda)=-\lambda\mu_{max}$
and $h_N(\lambda)=-\lambda\mu_{min}$ include all the lines $h(\lambda)=-\lambda\mu$.
Therefore, (\ref{eq:5_mar}) yields the two optima $\tilde\omega_1(\mmax)$ and $\tilde\omega_1(\mmin)$. However, these values must be equal as there is only one optimum, hence
\begin{equation}\label{eq:5_s05abbxx}
\displaystyle\frac{4\hat\omega_2\mu_{\max}}{(1+\hat\omega_2\mmax)^2}=\displaystyle\frac{4\hat\omega_2\mu_{\min}}{(1+\hat\omega_2\mmin)^2}
\end{equation}
or
\begin{equation}\label{eq:5_s05bbbxx}
\hat\omega_2=\displaystyle\frac{1}{\sqrt{\mmin\mmax}}
\end{equation}
which, because of (\ref{eq:5_s01bbxxaxax}), yields the optimum value for $\omega_2$ given by the second part of (\ref{eq:3u2ux}). Substituting the value of $\hat\omega_2$ in the expressions $\tilde\omega_1(\mmax)$ or $\tilde\omega_1(\mmin)$, given by the first or second part of the equality (\ref{eq:5_s05abbxx}), respectively, we obtain the optimum value of $\omega_1$ given by the first part of (\ref{eq:3u2ux}).
The spectral radius is given by
\begin{equation}\label{eq:5_s02bb0xx}
\rho{(\cl(\omega_1,\omega_2,a))}=\max \{|\tilde\lambda_1|,|\tilde\lambda_N|\}
\end{equation}
where $\tilde\lambda_1, \tilde\lambda_N$ are the abscissas of the
points of tangent of $h_1(\lambda), h_N(\lambda)$, respectively.
For the minimization of $\rho(\cl(\omega_1, \omega_2,a))$ with respect to $\omega_1$ we require
\begin{equation*}
|\tilde\lambda_1|=|\tilde\lambda_N|
\end{equation*}
or
\begin{equation}\label{eq:5_s2bbbxx}
\tilde\lambda_1=-\tilde\lambda_N=\left(1-\tilde\omega_1(\mmax)\right)^{1/2},
\end{equation}
where the last equality holds by the fact that $\tilde\lambda_1,\; \tilde\lambda_N$ are the abscissas of the tangents $h_1(\lambda)$ and $h_N(\lambda)$, respectively.
From (\ref{eq:5_s2bbbxxa}) and (\ref{eq:5_s2bbbxx}) it follows that
\begin{equation*}
\rho(\cl(\omega_1, \omega_2,a))=(1-\omega_{1_{opt}})^{1/2}
\end{equation*}
which, because of (\ref{eq:3u2ux}), yields (\ref{eq:4u2ux}). \qed
Theorem \ref{theo:opt_gsor} finds the optimum values of the relaxation parameters $\omega_1$ and $\omega_2$ of the GSOR($a$) method.
Letting $a=0$ in (\ref{eq:3u2ux}) we obtain the optima found also in \cite{BaiParlWang2005}.
Note that the parameter $a$ has no impact on the spectral radius of the GSOR($a$) method as one might have expected.
The algebraic approach in \cite{BaiParlWang2005} is similar to the one followed by \cite{Young} for determining the optimum of the sole parameter in the SOR method. In case of GSOR($a$), which has two parameters, there is an alternative less tedious algebraic approach (see \cite{Young} pp. 279-281). However, it remains to be verified whether either approach can be used to solve the problem of determining the optimum values of more than two parameters as is the case for the GMESOR($a$) method. Our approach follows the geometric approach of Varga \cite{Var1} for the determination of the optimum value of the parameter $\omega$ in SOR. It should be noted that this approach is also mentioned in \cite{Tayl} but without a proof.

\begin{cor}\label{cor:opt_GBSOR}
Consider the GBSOR($a$) method. Under the hypothesis of Theorem \ref{theo:opt_gsor} the spectral radius of the GBSOR($a$) method, $\rho(\cm(\omega_1, \omega_2, a))$, is minimized for any $a\neq 1+\sqrt{\mmin\mmax}$ at
\begin{equation}\label{eq:twa3u2u}
\omega_{1_{opt}}=\frac{4\sqrt{\mu_{min}\mu_{max}}}{(\sqrt{\mu_{min}}+\sqrt{\mu_{max}})^2}\ \ \mbox{and} \ \ \omega_{2_{opt}}=\frac{1}{(1-a)+\sqrt{\mu_{min}\mu_{max}}}
\end{equation}
and its corresponding value is
\begin{equation}\label{eq:twa4u2u}
\rho(\cm{(\omega_{1_{opt}},\omega_{2_{opt}},a)})=\frac{\sqrt{\mu_{max}}-\sqrt{\mu_{min}}}{\sqrt{\mu_{max}}+\sqrt{\mu_{min}}}.
\end{equation}
\end{cor}
{\bf{Proof}} We remark that the functional relationship (\ref{eq:20b}) of GBSOR($a$)
\noindent is the same as that of the GSOR($a$) method (\ref{eq:17}) with the only difference that now we have $1-a$ instead of $a$. Therefore, we have the same results as in Theorem \ref{theo:opt_gsor}, if we simply replace $a$ with $1-a$.\qed
If the matrix Q is symmetric negative definite, the optimum parameters and the minimum spectral radius for the GSOR($a$) method are given by the following theorem.
\begin{theo}\label{theo:opt_gsornon}
Consider the GSOR($a$) method. Let $A\in\mathbb{R}^{m\times m}$  be
symmetric positive definite, $B\in\mathbb{R}^{m\times n}$ be of
full column rank and $Q\in\mathbb{R}^{n\times n}$ be
symmetric negative definite. Denote the minimum and the maximum
eigenvalues of the matrix $J=Q^{-1}B^TA^{-1}B$ by $\mu_{min}$ and $\mu_{max}$, respectively.
Then the spectral radius of the GSOR($a$) method $\rho(\cl(\omega_1, \omega_2,a))$, when the matrix $J$ has negative eigenvalues, is minimized for any $a\neq\sqrt{\mmin\mmax}$ at
\begin{equation}\label{eq:3u2uxn}
\omega_{1_{opt}}=\frac{4\sqrt{\mu_{min}\mu_{max}}}{(\sqrt{|\mu_{min}|}+\sqrt{|\mu_{max}|})^2}\ \ \mbox{and} \ \ \omega_{2_{opt}}=\frac{1}{a-\sqrt{\mu_{min}\mu_{max}}}
\end{equation}
and its corresponding value is
\begin{equation}\label{eq:4u2uxnP}
\rho(\cl{(\omega_{1_{opt}},\omega_{2_{opt}},a)})=\frac{\sqrt{|\mu_{min}|}-\sqrt{|\mu_{max}|}}{\sqrt{|\mu_{min}|}+\sqrt{|\mu_{max}|}}.
\end{equation}
\end{theo}
{\bf{Proof}}{ { In this case $\mu < 0 $. Following a similar approach as in Theorem \ref{theo:opt_gsor}, we obtain (\ref{eq:3u2uxn}) and (\ref{eq:4u2uxnP}).}}   \qed

Under the hypothesis of Theorem \ref{theo:opt_gsornon} and if $a=0$, these results were also obtained in \cite{BaiParlWang2005}.
\begin{theo}\label{theo:opt_gbsornon}
Consider the GBSOR($a$) method. Let $A\in\mathbb{R}^{m\times m}$  be
symmetric positive definite, $B\in\mathbb{R}^{m\times n}$ be of
full column rank and $Q\in\mathbb{R}^{n\times n}$ be
symmetric negative definite. Denote the minimum and the maximum
eigenvalues of the matrix $J=Q^{-1}B^TA^{-1}B$ by $\mu_{min}$ and $\mu_{max}$, respectively.
Then the spectral radius of the GBSOR($a$) method $\rho(\cm(\omega_1, \omega_2,a))$, when the matrix $J$ has negative eigenvalues, is minimized for any $a\neq 1-\sqrt{\mmin\mmax}$ at
\begin{equation}\label{eq:3u2uxnb}
\omega_{1_{opt}}=\frac{4\sqrt{\mu_{min}\mu_{max}}}{(\sqrt{|\mu_{min}|}+\sqrt{|\mu_{max}|})^2}\ \ \mbox{and} \ \ \omega_{2_{opt}}=\frac{1}{1-a-\sqrt{\mu_{min}\mu_{max}}}
\end{equation}
and its corresponding value is
\begin{equation}\label{eq:4u2uxnPb}
\rho(\cl{(\omega_{1_{opt}},\omega_{2_{opt}},a)})=\frac{\sqrt{|\mu_{min}|}-\sqrt{|\mu_{max}|}}{\sqrt{|\mu_{min}|}+\sqrt{|\mu_{max}|}}.
\end{equation}
\end{theo}
{\bf{Proof}} We remark that the functional relationship (\ref{eq:20b}) of GBSOR($a$)
is the same as that of the GSOR($a$) method (\ref{eq:17}) with the only difference that now we have $1-a$ instead of $a$. Therefore, we can apply the results of Theorem \ref{theo:opt_gsornon} by replacing $a$ with $1-a$.\qed

}}

\subsubsection{The GMESOR($a$) method}
In the sequel we determine the optimum parameters for the GMESOR($a$) method.
\begin{theo}\label{theo:opt_1}
Consider the GMESOR($a$) method. Let $A\in\mathbb{R}^{m\times m}$ and $Q\in\mathbb{R}^{n\times n}$ be
symmetric positive definite and $B\in\mathbb{R}^{m\times n}$ be of
full column rank. Denote the minimum and the maximum eigenvalues of the matrix $J=Q^{-1}B^TA^{-1}B$ by
$\mu_{min}$ and $\mu_{max}$, respectively.
Then the spectral radius of the GMESOR($a$) method, $\rho(\ch(\tau_1, \tau_2,\omega_2,a))$, is minimized for any $a\neq -\sqrt{\mmin\mmax}$ at
\begin{equation}\label{eq:optgmesor1}
\omega_{2_{opt}}=\tau_{2_{opt}},
\end{equation}
\begin{equation}\label{eq:optgmesor2}
\tau_{1_{opt}}=\frac{4\sqrt{\mu_{min}\mu_{max}}}{(\sqrt{\mu_{min}}+\sqrt{\mu_{max}})^2}\ \ \mbox{and} \ \ \tau_{2_{opt}}=\frac{1}{a+\sqrt{\mu_{min}\mu_{max}}}
\end{equation}
and its corresponding value is
\begin{equation}\label{eq:optgmesor3}
\rho(\ch{(\tau_{1_{opt}},\tau_{2_{opt}},\omega_{2_{opt}},a}))=\frac{\sqrt{\mu_{max}}-\sqrt{\mu_{min}}}{\sqrt{\mu_{max}}+\sqrt{\mu_{min}}}.
\end{equation}
\end{theo}
{\bf{Proof}}
The functional relationship of the GMESOR($a$) method is given by (\ref{eq:10})
\noindent or
\begin{equation}\label{eq:5_s01bb}
(\lambda+\tau_1-1)(\lambda-1)=\frac{\tau_1(\omega_2-\tau_2-\lambda\omega_2)\mu}{1-a\omega_2}.
\end{equation}
The optimum values of $\tau_1,\;\tau_2\;\;\mbox{and}\;\; \omega_2$
will be determined such that
\begin{equation}\label{eq:5_s02bb}
\rho{(\ch(\tau_1,\tau_2,\omega_2, a))}=\max_{\mu_{min}\leq\mu\leq\mu_{max}}|\lambda|
\end{equation}
is minimum.
The real roots of (\ref{eq:10}) are the intersection points
of the parabola
\begin{equation}\label{eq:5_s03bb}
g(\lambda)=\frac{(\lambda+\tau_1-1)(\lambda-1)}{\tau_1}
\end{equation}
and the straight lines
\begin{equation}\label{eq:5_s04bb}
h(\lambda)=\frac{\omega_2-\tau_2-\lambda\omega_2}{1-a\omega_2}\mu, \;\; 0<\mu_{min}\leq\mu\leq\mu_{max}.
\end{equation}
Following a similar argument as in \cite{Var1} page 111,
$h(\lambda)$ are straight lines through the point $\left(0,\displaystyle\frac{\omega_2-\tau_2}{1-a\omega_2}\mu\right)$
and $g_{\tau_1}(\lambda)$ is a parabola passing through the points (1,0) and $(1-\tau_1,0)$ (see figure 3). 
\begin{figure}[ht]
  \begin{center}
   \scalebox{0.6}{\includegraphics{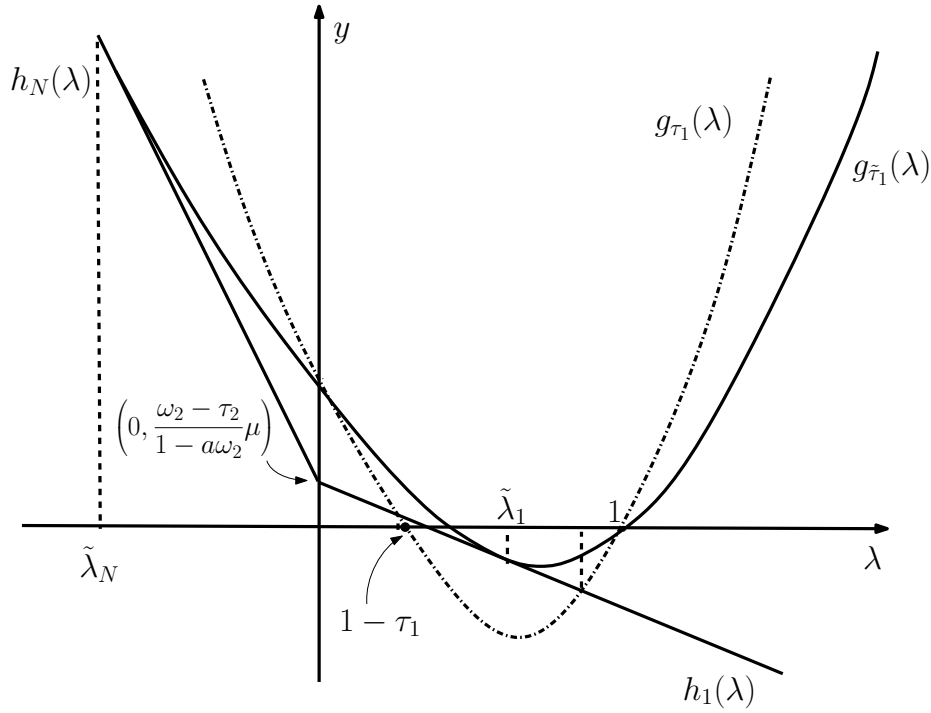}}
  \end{center}\label{fig:5-02}
  \caption{Graphs of $g_{\tau_1}(\lambda), h_1(\lambda)$ and $h_N(\lambda)$ in case the roots of (\ref{eq:5_s01bb}) are real.}
\end{figure}
The spectral radius is given by
\begin{equation}\label{eq:5_s02bb0}
\rho{(\ch(\tau_1,\tau_2,\omega_2,a))}=\max \{|\tilde\lambda_1|,|\tilde\lambda_N|\}
\end{equation}
where $\tilde\lambda_1, \tilde\lambda_N$ are the abscissas of the
points of tangent of $h_1(\lambda), h_N(\lambda)$, respectively, where now $h_1(\lambda)=(\omega_2-\tau_2-\lambda\omega_2)\mmax$ and $h_N(\lambda)=(\omega_2-\tau_2-\lambda\omega_2)\mmin$.
Therefore,
\begin{equation}\label{eq:5_s02bb1}
|\tilde\lambda_1|=\left(1-\tau_1+\tau_1\frac{\tau_2-\omega_2}{1-a\omega_2}\mu_{max}\right)^{1/2}
\end{equation}
\noindent and
\begin{equation}\label{eq:5_s02bb2}
|\tilde\lambda_N|=\left(1-\tau_1+\tau_1\frac{\tau_2-\omega_2}{1-a\omega_2}\mu_{min}\right)^{1/2}.
\end{equation}
From (\ref{eq:5_s02bb0}) it follows that the minimum value of
$\rho{(\ch(\tau_1,\tau_2,\omega_2, a))}$ is attained when
\begin{equation}\label{eq:5_s06bbaaa}
|\tilde\lambda_1|=|\tilde\lambda_N|
\end{equation}
which, because of (\ref{eq:5_s02bb1}) and (\ref{eq:5_s02bb2}), implies

\begin{equation}\label{eq:s06bb}
\omega_2=\tau_2.
\end{equation}
In case $\tilde\lambda_1$ and $\tilde\lambda_N$ are the two conjugate complex roots of (\ref{eq:5_s01bb}), it follows that (\ref{eq:5_s06bbaaa}) must also hold for $\rho{(\ch(\tau_1,\tau_2,\omega_2,a))}$ to be minimized.
So, (\ref{eq:s06bb}) holds if either (\ref{eq:5_s01bb}) has real or conjugate complex roots.
However, if (\ref{eq:s06bb}) holds, then (\ref{eq:10}) becomes
\begin{equation*}
\lambda^2+\lambda\left(\tau_1-2+\tau_1\hat\tau_2\mu\right)+1-\tau_1=0,
\end{equation*}
\noindent which is the functional relationship of the GSOR with
\begin{equation}\label{eq:5_s06bb5x1}
\hat\tau_2=\frac{\tau_2}{1-a\tau_2}.
\end{equation}
Therefore the optimum values of $\tau_1$ and $\hat\tau_2$ are given by $\omega_{opt}$ and $\tau_{opt}$ of \cite{BaiParlWang2005}, respectively, whereas the minimum value of $\rho{(\ch(\tau_1,\tau_2,\omega_2, a))}$ is given by $\rho{(\ch(\omega_{opt},\tau_{opt}))}$ of \cite{BaiParlWang2005}. Finally, using (\ref{eq:5_s06bb5x1}) we find (\ref{eq:optgmesor2}). \qed
\noindent
So, for the optimum values of its parameters, GMESOR($a$) degenerates to the { {GSOR($a$)}} method.

{ {
\begin{cor}\label{cor:opt_gesor}
Consider the GESOR($a$) method. Let $A\in\mathbb{R}^{m\times m}$ and $Q\in\mathbb{R}^{n\times n}$ be
symmetric positive definite and $B\in\mathbb{R}^{m\times n}$ be of
full column rank. Denote the minimum and the maximum eigenvalues of the matrix $J=Q^{-1}B^TA^{-1}B$ by
$\mu_{min}$ and $\mu_{max}$, respectively. Then the spectral radius of the GESOR($a$) method, $\rho(\ch(\tau, \omega_2, a))$, is minimized at
\begin{equation}\label{eq:tw}
\omega_{2_{opt}}=\tau_{opt}
\end{equation}
and
\begin{equation}\label{eq:twa2}
\tau_{{opt}}=\frac{4\sqrt{\mu_{min}\mu_{max}}}{(\sqrt{\mu_{min}}+\sqrt{\mu_{max}})^2},\;\; a_{opt}=\frac{1}{\tau_{opt}}-\sqrt{\mmax\mmin}
\end{equation} 
and its corresponding value is
\begin{equation}\label{eq:twa4}
\rho(\ch{(\tau_{opt},\omega_{2_{opt}}},a_{opt}))=\frac{\sqrt{\mu_{max}}-\sqrt{\mu_{min}}}{\sqrt{\mu_{max}}+\sqrt{\mu_{min}}}.
\end{equation}
\end{cor}
{\bf{Proof}} Recall that GESOR($a$) is obtained by setting $\tau_1=\tau_2$ in GMESOR($a$). Therefore, (\ref{eq:tw}) and (\ref{eq:twa2}) are obtained by (\ref{eq:optgmesor1}) and (\ref{eq:optgmesor2}), respectively, where now we require $\tau_{opt}=\tau_{1_{opt}}=\tau_{2_{opt}}$.\qed

\begin{cor}\label{cor:opt_GMEBSOR}
Consider the GMEBSOR($a$) method. Under the hypothesis of Theorem \ref{theo:opt_1} the spectral radius of the GMEBSOR($a$) method, $\rho(\ck(\tau_1, \tau_2, \omega_1, \omega_2, a))$, is minimized at
\begin{equation}\label{eq:twu}
\omega_{1_{opt}}=\tau_{1_{opt}},
\end{equation}
where
\begin{equation}\label{eq:twa3u}
\tau_{1_{opt}}=\frac{4\sqrt{\mu_{min}\mu_{max}}}{(\sqrt{\mu_{min}}+\sqrt{\mu_{max}})^2}\ \ \mbox{and} \ \ \tau_{2_{opt}}=\frac{1-(1-a)\omega_2}{\sqrt{\mu_{min}\mu_{max}}}
\end{equation}
and its corresponding value is
\begin{equation}\label{eq:twa4u}
\rho(\ck{(\tau_{1_{opt}},\tau_{2_{opt}},\omega_{1_{opt}},\omega_2, a}))=\frac{\sqrt{\mu_{max}}-\sqrt{\mu_{min}}}{\sqrt{\mu_{max}}+\sqrt{\mu_{min}}}.
\end{equation}
\end{cor}
{\bf{Proof}} Following a similar approach as in Theorem \ref{theo:opt_1}, using the functional relationship (\ref{eq:18}) and requiring $|\tilde\lambda_1|=|\tilde\lambda_N|$ we find
\begin{equation}\label{eq:optGMEBSOR1}
\omega_1=\tau_1.
\end{equation}
Therefore, (\ref{eq:18}) because of (\ref{eq:optGMEBSOR1}) becomes
\begin{equation}\label{eq:optGMEBSOR2}
\lambda^2+\lambda\left(\tau_1-2+\tau_1\hat\tau_2\mu\right)+1-\tau_1=0
\end{equation}
\noindent
with
\begin{equation}\label{eq:optGMEBSOR2x}
\hat\tau_2=\frac{\tau_2}{1-(1-a)\omega_2},
\end{equation}
which is the functional relationship of the GSOR($a$) method (see (\ref{eq:17}))
\noindent with the only difference that now we have $1-a$ instead of $a$ in $\hat\tau_2$, hence (\ref{eq:twa3u}) and (\ref{eq:twa4u}) hold because of Theorem \ref{theo:opt_gsor}. \qed
Note that although the GMEBSOR($a$) method has four parameters instead of three as in the GMESOR($a$) method, both methods have the same minimum spectral radius.

\begin{cor}\label{cor:opt_GEBSOR}
Consider the GEBSOR($a$) method. Under the hypothesis of Theorem \ref{theo:opt_1} the spectral radius of the GEBSOR($a$) method, $\rho(\ck(\tau, \omega_2, a))$, is minimized at
\begin{equation}\label{eq:twuu1}
\omega_{1_{opt}}=\tau_{opt}
\end{equation}
and
\begin{equation}\label{eq:twau2}
\tau_{{opt}}=\frac{4\sqrt{\mu_{min}\mu_{max}}}{(\sqrt{\mu_{min}}+\sqrt{\mu_{max}})^2},\;\; \omega_{2_{opt}}=\frac{1-\tau_{opt}\sqrt{\mmin\mmax}}{1-a}
\end{equation}
and its corresponding value is
\begin{equation}\label{eq:twau4}
\rho(\ck{(\tau_{opt},\omega_{2_{opt}}, a)})=\frac{\sqrt{\mu_{max}}-\sqrt{\mu_{min}}}{\sqrt{\mu_{max}}+\sqrt{\mu_{min}}}.
\end{equation}
\end{cor}
{\bf{Proof}} Recall that GEBSOR($a$) is obtained by setting $\tau_1=\tau_2$ in GMEBSOR($a$). Therefore, (\ref{eq:twuu1}) and (\ref{eq:twau2}) are obtained by (\ref{eq:twu}) and (\ref{eq:twa3u}), respectively, where now we require $\tau_{opt}=\tau_{1_{opt}}=\tau_{2_{opt}}$.\qed

Our analysis so far shows that all the studied iterative methods (GMESOR($a$), GMEBSOR($a$)) have also the same rate of convergence as the PCG method for the optimum values of their parameters (see Theorems \ref{theo:opt_gsor}, \ref{theo:opt_1} and corollary \ref{cor:opt_GMEBSOR}).
}}

\section{The Generalized Modified Preconditioned Simultaneous Displacement (GMPSD) method}
The Preconditioned Simultaneous Displacement (PSD) method was introduced in \cite{EvMis1}.  When the coefficient matrix A is two-cyclic
the Modified PSD (MPSD) method was studied in \cite{LoukMis2}, \cite{MisEv2}. Motivated by our previous work we introduce the Generalized Modified PSD (GMPSD) method and study its convergence rate for the numerical solution of the augmented linear system (\ref{eq:01})-(\ref{eq:02}).

\subsection{The functional relationship}
In the sequel, we let the preconditioning matrix $R$ be the product of the lower triangular part with the upper triangular part of $\ca$ in an attempt to obtain a better approximation of $\ca$ and consequently an increase in the rate of convergence of the corresponding iterative method. Let
\begin{equation}\label{eq:21}
R=(\cd-\Omega\cl)\cd^{-1}(\cd-\Omega\cu).
\end{equation}
From (\ref{eq:06}) and (\ref{eq:21}) it follows that the iteration matrix of (\ref{eq:05}) now is
\begin{equation}\label{eq:22}
\cg{(\tau_1,\tau_2,\omega_1,\omega_{2},a)}=I-(\cd-\Omega \cu)^{-1}D(\cd-\Omega \cl)^{-1}T\ca
\end{equation}
whereas $\eta(\tau_1, \tau_2)$ in (\ref{eq:06}) corresponds to
\begin{equation}\label{eq:23}
\gamma{(\tau_1,\tau_2,\omega_1,\omega_{2},a)}=(\cd-\Omega \cu)^{-1}D(\cd-\Omega \cl)^{-1}Tb.
\end{equation}
Note that this method  has four parameters $\tau_1, \tau_2, \omega_1\; \mbox{and}\; \omega_2$ instead of three in the GMESOR method.
The iterative scheme given by (\ref{eq:05}), (\ref{eq:22}) and (\ref{eq:23})  will be referred to as the Generalized Modified Preconditioned Simultaneous Displacement (GMPSD) method. For $(\cd-\Omega \cu)^{-1}D(\cd-\Omega \cl)^{-1}$ to exist we require
\begin{equation}\label{eq:5_06bbbblu}
\det[(\cd-\Omega\cl)\cd^{-1}(\cd-\Omega\cu)]\neq 0.
\end{equation}
Because of (\ref{eq:04})
\begin{equation}\label{eq:5_06bbbclu}
R=(\cd-\Omega\cl)\cd^{-1}(\cd-\Omega\cu)=\left(
\begin{array}{cc}
A & \omega_1 B\\
-\omega_2 B^T & (1-a\omega_2)[1-(1-a)\omega_2] Q{{-\omega_1\omega_2B^TA^{-1}B}}
\end{array}\right).
\end{equation}
Therefore,
\begin{equation*}
\det(\cd-\Omega\cl)\cd^{-1}(\cd-\Omega\cu)=(1-a\omega_2)^{n}[1-(1-a)\omega_2]^{n}\det{(A)}\det{(Q)}\neq 0
\end{equation*}
or
\begin{equation}\label{eq:5_100}
a\neq \frac12 \;\;\mbox{and}\;\; \omega_2\neq 2
\end{equation}
since the matrix $A$ is symmetric positive definite and the matrix $Q$ is nonsingular.
The GMPSD method has the following algorithmic form.
\\ \noindent
{\sc The GMPSD Method:}{\em{ Let $Q\in\mathbb{R}^{n\times n}$ be a nonsingular and symmetric matrix.
Given initial vectors $x^{(0)}\in\mathbb{R}^m$ and
$y^{(0)}\in\mathbb{R}^n$, and relaxation factors $\tau_1,\ \tau_2\neq 0,  \
\omega_1, \omega_2, a\in\mathbb{R}$ with $a\neq \frac12 \;\; \mbox{and}\;\; \omega_2\neq 2$.
For $k=0,1,2,...$ until the
iteration sequence
$\{({x^{(k)}}^T, {y^{(k)}}^T)^T\}$ is convergent, compute
\\ \newline
\hspace{0.5cm}${y^{(k+1)}=y^{(k)}+
\frac{1}{(1-a\omega_2)[1-(1-a)\omega_2]}Q^{-1}\left\{B^T[(\tau_2-\tau_1\omega_2)x^{(k)}+\tau_1\omega_2A^{-1}(b_{{1}}-By^{(k)})]{{-}}\tau_2b_2\right\}}$
\normalsize
\\
\hspace{0.3cm}$    {x^{(k+1)}=(1-\tau_1)x^{(k)}+ A^{-1}\left\{B\left[(\omega_1-\tau_1)y^{(k)}-\omega_1y^{(k+1)}\right]+\tau_1b_1\right\}},
$ \\ \newline
\noindent where Q is an  approximation of the Schur complement matrix $B^TA^{-1}B$.}}
\\\\
Note that in the above algorithm we first compute $y^{(k+1)}$ and then $x^{(k+1)}$, whereas in the GMESOR method we had the reverse computations. { {If $\tau=\tau_1=\tau_2$ and $\omega=\omega_1=\omega_2$ we have the GPSD method. }}



\noindent If $\omega_2=0$ then the algorithmic form of the GMPSD method simplifies to
\begin{equation}\label{GMPSD0}
\begin{array}{ll}
y^{(k+1)}=y^{(k)}+\tau_2 Q^{-1}(B^T x^{(k)}{{-}}b_2)&\\
x^{(k+1)}=(1-\tau_1)x^{(k)}+\tau_1 A^{-1}(b_1-By^{(k+1)})&
\end{array}
\end{equation}
The above form is the same as that of the GSOR method if we use $\cd-\Omega \cu$ instead of $\cd-\Omega \cl$ as the preconditioned matrix in the GSOR method and will be referred as the simplified GMPSD method.
\noindent In the following theorem we find the functional relationship for the GMPSD method between the eigenvalues $\lambda$ of the iteration matrix $\cg(\tau_1, \tau_2, \omega_1, \omega_2, a)$  and the eigenvalues $\mu$ of the matrix $J$.
\begin{theo}\label{theo:func_rel2}
Let $A\in\mathbb{R}^{m\times m}$ be symmetric positive definite, $B\in\mathbb{R}^{m\times n}$ be of full column rank and $Q\in\mathbb{R}^{n\times n}$ be
nonsingular and symmetric. If $\lambda\neq 1-\tau_1$ is an eigenvalue of the matrix
$\cg(\tau_1, \tau_2, \omega_1, \omega_2, a)$ and if $\mu$ satisfies
\footnotesize{
\begin{equation}
\lambda^2+\lambda\left(\tau_1-2+\displaystyle\frac{\tau_1\omega_2+\tau_2\omega_1-\tau_1\omega_1\omega_2}{(1-a\omega_2)[1-(1-a)\omega_2]}\mu\right)+1-\tau_1+\displaystyle\frac{\tau_1\tau_2-\tau_1\omega_2-\tau_2\omega_1+\tau_1\omega_1\omega_2}{(1-a\omega_2)[1-(1-a)\omega_2]}\mu=0, \label{eq:24}
\end{equation}
}\normalsize
\noindent where $a\neq \frac12 \; \mbox{and}\; \omega_2\neq 2$, then $\mu$ is an eigenvalue of the matrix $J=Q^{-1}B^TA^{-1}B$.
Conversely, if $\mu$ is an eigenvalue of $J$ and if
$\lambda\neq1-\tau_1$ satisfies (\ref{eq:24}), then $\lambda$ is an
eigenvalue of $\cg(\tau_1, \tau_2, \omega_1, \omega_2, a)$. In addition, $\lambda=1-\tau_1$ is an
eigenvalue of $\cg(\tau_1, \tau_2, \omega_1, \omega_2, a)$ (if $m>n$) with the corresponding
eigenvector $(x^T,0)^T$, where $x\in\cn(B^T)$.
\end{theo}
{\bf{Proof}}
{{Clearly, the eigenvalues $\mu$ of the matrix $J=Q^{-1}B^TA^{-1}B$ are real and non-zero. Let $\lambda$ be a nonzero eigenvalue of the iteration matrix $\cg{(\tau_1,\tau_2,\omega_1,\omega_{2},a)}$ defined in (\ref{eq:22}), and $(x,y)^T\in\mathbb{R}^{m+n}$ be the corresponding eigenvector. Then, we have that
\begin{equation*}
\cg{(\tau_1,\tau_2,\omega_1,\omega_{2},a)}\left(
\begin{array}{c}
x\\y
\end{array}
\right)=\lambda\left(
\begin{array}{c}
x\\y
\end{array}
\right)
\end{equation*}
or because of (\ref{eq:22})
\begin{equation}\label{eq:11}
[(\cd-\Omega\cl)\cd^{-1}(\cd-\Omega\cu)-T\ca]\left(
\begin{array}{c}
x\\y
\end{array}
\right)=\lambda (\cd-\Omega\cl)\cd^{-1}(\cd-\Omega\cu) \left(
\begin{array}{c}
x\\y
\end{array}
\right).
\end{equation}
\\ \noindent From (\ref{eq:11}), because of (\ref{eq:04}), we have that
\begin{eqnarray*}
\left(
\begin{array}{ccc}
(1-\tau_1)A&&(\omega_1-\tau_1) B
\\
(\tau_2-\omega_2)B^T&&(1-a\omega_2)[1-(1-a)\omega_2]Q-\omega_1\omega_2B^TA^{-1}B
\end{array}
\right)\left(
\begin{array}{c}
x\\y
\end{array}
\right)\hspace{1cm}\\
=\lambda \left(
\begin{array}{ccc}
A&&\omega_1 B
\\
-\omega_2B^T&&(1-a\omega_2)[1-(1-a)\omega_2]Q-\omega_1\omega_2B^TA^{-1}B
\end{array}
\right) \left(
\begin{array}{c}
x\\y
\end{array}
\right).
\end{eqnarray*}

\noindent Decoupling we have that
\begin{equation*}
\begin{cases}
(1-\tau_1)Ax+(\omega_1-\tau_1) By=\lambda Ax+\lambda\omega_1 By\\
(\tau_2-\omega_2)B^Tx+\{(1-a\omega_2)[1-(1-a)\omega_2]Q-\omega_1\omega_2B^TA^{-1}B\}y\\\hspace{1cm}=-\lambda \omega_2 B^T x+\lambda\{(1-a\omega_2)[1-(1-a)\omega_2]Q-\omega_1\omega_2B^TA^{-1}B\}y
\end{cases}
\end{equation*}
or equivalently
\begin{equation}\label{eq:12}
\begin{cases}
(1-\tau_1-\lambda)x=[(\lambda-1)\omega_1+\tau_1] A^{-1}By\\
(\tau_2-\omega_2+\lambda\omega_2)Q^{-1}B^Tx=(\lambda-1)\left\{(1-a\omega_2)[1-(1-a)\omega_2]{{I}}-\omega_1\omega_2J\right\}y.
\end{cases}
\end{equation}
From the first equality in (\ref{eq:12}) we get
\begin{equation*}
(1-\tau_1-\lambda)Q^{-1}B^Tx=[(\lambda-1)\omega_1+\tau_1]Jy,
\end{equation*}
and hence, when $\lambda\neq 1-\tau_1$,
\begin{equation}\label{eq:13}
Q^{-1}B^Tx=\frac{(\lambda-1)\omega_1+\tau_1}{1-\tau_1-\lambda} Jy.
\end{equation}
It then follows from (\ref{eq:13}) and the second equality in (\ref{eq:12}) that
\begin{eqnarray*}
(\lambda-1)(1-a\omega_2)[1-(1-a)\omega_2](1-\tau_1-\lambda)y\hspace{4cm}\\ =\{[(\lambda-1)\omega_2+\tau_2][(\lambda-1)\omega_1+\tau_1] +(\lambda-1)(1-\tau_1-\lambda)\omega_1\omega_2\}J y.
\end{eqnarray*}
\noindent If $\lambda=1-\tau_1\neq0$, then from the first and the second equality of (\ref{eq:12}) we have, respectively, $By=0$ and $\tau_1\{(1-a\omega_2)[1-(1-a)\omega_2]Q-\omega_1\omega_2B^TA^{-1}B\}y=(\tau_1\omega_2-\tau_2)B^T x$. It then follows that $y=0$ and $x\in\cn(B^T)$, where $\cn(B^T)$ is the null space of the matrix $B^T$. Hence, {$\lambda=1-\tau_1$} is an eigenvalue of $\cg{(\tau_1,\tau_2,\omega_1,\omega_{2},a)}$ with the corresponding eigenvector $(x^T, 0)^T$, where $x\in\cn(B^T)$. Therefore, the eigenvalues $\lambda$ (except for $\lambda=1-\tau_1$) of the matrix $\cg{(\tau_1,\tau_2,\omega_1,\omega_{2},a)}$ and the eigenvalues $\mu$ of the matrix $J$ satisfy the functional relationship
\begin{eqnarray*}
(\lambda-1)(1-a\omega_2)[1-(1-a)\omega_2](1-\tau_1-\lambda)\hspace{4cm}\\ =\left\{[(\lambda-1)\omega_2+\tau_2][(\lambda-1)\omega_1+\tau_1] +(\lambda-1)(1-\tau_1-\lambda)\omega_1\omega_2\right\}\mu.
\end{eqnarray*}
This means that $\lambda$ satisfies the quadratic equation (\ref{eq:24}). \qed }}

{ {
\begin{cor} Let $A\in\mathbb{R}^{m\times m}$ be symmetric positive definite, $B\in\mathbb{R}^{m\times n}$ be of full column rank and $Q\in\mathbb{R}^{n\times n}$ be
nonsingular and symmetric.
\ \newline \noindent
1. The nonzero eigenvalues of the iteration matrix $\cg(\tau, \omega_1, \omega_2, a)$ of the
GMPSD(3) method are given by $\lambda=1-\tau$ or if $a\neq \frac12 \;\; \mbox{and}\;\; \omega_2\neq 2$ by
\begin{equation}
\lambda^2+\lambda\left({\tau-2}+\frac{\tau\hat\omega}{(1-a\omega_2)[1-(1-a)\omega_2]}\mu\right) +1-\tau+\frac{\tau(\tau-\hat\omega)}{(1-a\omega_2)[1-(1-a)\omega_2]}\mu=0
\label{eq:26}
\end{equation}
where
\begin{equation}\label{eq:26b}
\hat\omega=\omega_1+\omega_2-\omega_1\omega_2.
\end{equation}
2. The nonzero eigenvalues of the iteration matrix $\mathcal{S}(\omega_1, \omega_2, a)$ of the
GMSSOR method are given by $\lambda=1-\hat\omega$ or if $a\neq \frac12 \;\; \mbox{and}\;\; \omega_2\neq 2$ by
\begin{equation}
\lambda^2+\lambda\left({\hat\omega-2}+\frac{\hat\omega^2}{(1-a\omega_2)[1-(1-a)\omega_2]}\mu\right)+1-\hat\omega=0 \label{eq:27}
\end{equation}
where $\hat\omega$ is given by (\ref{eq:26b}).\\ \noindent
The nonzero eigenvalues of the iteration matrix $\cg(\tau, \omega, a)$ of the
GPSD method are given by $\lambda=1-\tau$ or if $a\neq \frac12 \;\; \mbox{and}\;\; \omega\neq 2$ by
\begin{equation}
\lambda^2+\lambda\left({\tau-2}+\frac{\tau\hat\omega}{(1-a\omega)[1-(1-a)\omega]}\mu\right) +1-\tau+\frac{\tau(\tau-\hat\omega)}{(1-a\omega)[1-(1-a)\omega]}\mu=0
\label{eq:28}
\end{equation}
where now
\begin{equation}\label{eq:28b}
\hat\omega=\omega(2-\omega).
\end{equation}
and
\begin{equation}\label{eq:28c}
 a\omega\neq 1 \;\;\mbox{and} \;\; (1-a)\omega\neq 1
\end{equation}
4. { {The nonzero eigenvalues of}} the iteration matrix $\mathcal{S}(\omega, a)$ of the
GSSOR method are given by $\lambda=1-\hat\omega$ or if $a\neq \frac12 \;\; \mbox{and}\;\; \omega\neq 2$ by
\begin{equation}
\lambda^2+\lambda\left({\hat\omega-2}+\frac{\hat\omega^2}{(1-a\omega)[1-(1-a)\omega]}\mu\right)+1-\hat\omega=0 \label{eq:29}
\end{equation}
where $\hat\omega$ is given by (\ref{eq:28b}).
\end{cor}

{\bf{Proof}} The iteration matrix $\cg(\tau,\omega_1,\omega_2, a)$ of the GMPSD(3) is obtained by letting $\tau=\tau_1=\tau_2$ in $\cg(\tau_1, \tau_2, \omega_1, \omega_2, a)$ given by (\ref{eq:22}).
Using the matrix $\cg(\tau_1, \tau_2, \omega_1, \omega_2, a)$ and following a similar approach as in the proof of Theorem \ref{theo:func_rel2} we find the functional relationship (\ref{eq:26}). Similarly, we find (\ref{eq:27}), (\ref{eq:28}) and (\ref{eq:29}).\qed

}}
\subsection{Convergence}
If the matrix $Q$ is positive definite and $a=0$ sufficient conditions for the GMPSD method to converge are given by the following theorem.
\begin{theo}\label{theo:conv1bb}
Consider the GMPSD method. Let $A\in\mathbb{R}^{m\times m}$ and $Q\in\mathbb{R}^{n\times n}$ be
symmetric positive definite and $B\in\mathbb{R}^{m\times n}$ be of
full column rank.
Denote the minimum and the maximum eigenvalues of the matrix $J=Q^{-1}B^TA^{-1}B$ by $\mu_{min}$ and $\mu_{max}$, respectively. Then,
$\rho({\cg}{(\tau_1,\tau_2,\omega_1,\omega_2)})<1$ if the
parameters $\tau_1,\tau_2,\omega_1\; \mbox{and}\; \omega_2$ lie in the region defined in the cases of Table \ref{tab10} with $0<\tau_1<2$ and
\begin{equation}\label{eq:conv_gmpsd}
\begin{array}{ll}
\omega_{11}^*(\mu)=\displaystyle\frac{\tau_1(2\omega_2-\tau_2)}{2(\tau_1\omega_2-\tau_2)}+\displaystyle\frac{(\tau_1-2)(1-\omega_2)}{\tau_1\omega_2-\tau_2}\displaystyle\frac{1}{\mu}, \;\;\;&\;\;\; \omega_{21}^*=\displaystyle\frac{\tau_2}{\tau_1},\\
 \omega_{12}^*(\mu)= \displaystyle\frac{\tau_1(\omega_2-\tau_2)}{\tau_1\omega_2-\tau_2}+\displaystyle\frac{\tau_1(1-\omega_2)}{\tau_1\omega_2-\tau_2}\displaystyle\frac{1}{\mu}, \;\;\;& \;\;\;\omega_{22}^*(\mu)=1-\displaystyle\frac{\tau_1\tau_2\mu}{4}.
\end{array}
\end{equation}
\begin{table}[htbp]
\caption{Sufficient conditions for the
 GMPSD method to converge. } \label{tab10}
\begin{tabular}{|c|c|c|c|c|c|c|c|c|c|c|}
\hline
   $\mbox{Cases}$&$\omega_2 - \mbox{Domain}$  &$\omega_1 - \mbox{Domain}$
   &${\tau_2 - \mbox{Domain} }$\\
    \hline\hline
1& $ \omega_{21}^*<\omega_2<\omega_{22}^*(\mu_{max}) $& $\omega_{11}^*(\mu_{max})<\omega_1<\omega_{12}^*(\mu_{min})$ &
 \\
 \cline{1-3} 2&$\omega_2<\omega_{21}^* $
 & &$0<\tau_2<{\displaystyle\frac{4\tau_1 }{4+\tau_1^2\mmax}
}$
 \\ \cline{1-2} \cline{4-4} 3&  $\omega_2<\omega_{22}^*(\mu_{max}) $
 & $\omega_{12}^*(\mu_{min})<\omega_1<\omega_{11}^*(\mu_{max})$ & ${\displaystyle\frac{4\tau_1}{4+\tau_1^2\mmin}
}<\tau_2$
 \\ \cline{1-2} \cline{4-4} 4&$1<\omega_2<\omega_{22}^*(\mmin)$
 & &$\tau_2<0$
 \\ \hline
\end{tabular}
\end{table}
\end{theo}
{\bf{Proof}} Recall that $\lambda=1-\tau_1\neq 0$ is an eigenvalue of $\cg{(\tau_1,\tau_2,\omega_1,\omega_{2})}$ and if $\lambda\neq 1-\tau_1$ then the eigenvalues of $\cg{(\tau_1,\tau_2,\omega_1,\omega_{2})}$ are given by (\ref{eq:24}). If $\lambda=1-\tau_1\neq 0$, then the GMPSD method is convergent if and only if $|\lambda|<1$, that is $|1-\tau_1|<1$, or
\begin{equation} \label{eq:45}
0<\tau_1<2.
\end{equation}
If $\lambda\neq 1-\tau_1$, then (\ref{eq:24}) holds and by Lemma 2.1 page 171 of \cite{Young}, it follows that the GMPSD method is convergent if and only if
\begin{equation}\label{eq:3g2}
|c|<1 \ \ \mbox{and}\ \ |b|<1+c
\end{equation}
where
\begin{equation} \label{eq:46}
c=1-\tau_1+\frac{\tau_1\omega_1\omega_2-\tau_1\omega_2-\tau_2\omega_1+\tau_1\tau_2}{1-\omega_2}\mu
\end{equation}
and
{{
\begin{equation} \label{eq:47}
b=1+c-\frac{\tau_1\tau_2}{1-\omega_2}\mu.
\end{equation}
}}
From the first inequality of (\ref{eq:3g2}) it follows that
\begin{equation}\label{eq:35}
0<1+c<2.
\end{equation}
From the second inequality of (\ref{eq:3g2}), because of (\ref{eq:47}), we have
%
\begin{equation} \label{eq:48}
0<\frac{\tau_1\tau_2\mu}{2(1-\omega_2)}<1+c.
\end{equation}
Combining (\ref{eq:35}) and (\ref{eq:48}) it follows that
\begin{equation} \label{eq:49}
0<\frac{\tau_1\tau_2\mu}{2(1-\omega_2)}<1+c<2.
\end{equation}
In order for (\ref{eq:49}) to hold we must have that
\begin{equation*} 
0<\frac{\tau_1\tau_2\mu}{2(1-\omega_2)}<2,
\end{equation*}
or because of (\ref{eq:45})
\begin{equation} \label{eq:50}
0<\frac{\tau_2}{1-\omega_2}<\frac{4}{\tau_1\mu}.
\end{equation}
Inequalities (\ref{eq:49}), because of (\ref{eq:46}), become
\begin{equation} \label{eq:51}
\frac{\tau_1(2\omega_2-\tau_2)\mu}{2(1-\omega_2)}+\tau_1-2<\omega_1\frac{\tau_1\omega_2-\tau_2}{1-\omega_2}\mu<\tau_1+\frac{\tau_1(\omega_2-\tau_2)\mu}{1-\omega_2}.
\end{equation}
In the sequel we distinguish the following two cases to study (\ref{eq:51}). Case I: $\tau_2>0$ and $1-\omega_2>0$ and Case II: $\tau_2<0$ and $1-\omega_2<0$. In addition, we distinguish the following two subcases for each of the above cases. (i): $\tau_1\omega_2-\tau_2>0$ and (ii): $\tau_1\omega_2-\tau_2<0$. Next, we will study only the subcase (i) of Case I, since the other cases can be treated similarly.
\noindent For this case, we have that
\begin{equation} \label{eq:52}
\frac{\tau_2}{\tau_1}<\omega_2<1, \;\;\; \mbox{if}\;\;\; 0<\tau_2<\tau_1
\end{equation}
and from the second part of (\ref{eq:50})
\begin{equation} \label{eq:53}
\omega_2<1-\frac{\tau_1\tau_2\mu}{4}.
\end{equation}
From (\ref{eq:52}) and (\ref{eq:53}) it follows that
\begin{equation*} 
\frac{\tau_2}{\tau_1}<\omega_2<\min{\left\{1, 1-\frac{\tau_1\tau_2\mu}{4}\right\}},\;\;\; 0<\tau_2<\tau_1
\end{equation*}
or
\begin{equation} \label{eq:54}
\frac{\tau_2}{\tau_1}<\omega_2<1-\frac{\tau_1\tau_2\mu}{4},\;\;\; 0<\tau_2<\tau_1
\end{equation}
which holds if $\frac{\tau_2}{\tau_1}<1-\frac{\tau_1\tau_2\mu}{4}$. Therefore, we have that (\ref{eq:54}) holds if
\begin{equation} \label{eq:55}
\omega_{21}^*<\omega_2<\omega_{22}^*(\mu), \;\;\; 0<\tau_2<\frac{4\tau_1}{4+\tau_1^2\mu}
\end{equation}
where $\omega_{21}^*, \omega_{22}^*(\mu)$ are given by (\ref{eq:conv_gmpsd}). Furthermore, from (\ref{eq:51}), we have that
\begin{equation} \label{eq:56}
\omega_{11}^*(\mu)<\omega_1<\omega_{12}^*(\mu)
\end{equation}
where $\omega_{11}^*(\mu), \omega_{12}^*(\mu)$ are given by (\ref{eq:conv_gmpsd}). Studying the monotonicity of $\omega_{22}^*(\mu), \omega_{11}^*(\mu)$ and $\omega_{12}^*(\mu)$ with respect to $\mu$ we have that
$sign\frac{\partial\omega_{22}^*(\mu)}{\partial\mu}=-1$, $sign\frac{\partial\omega_{11}^*(\mu)}{\partial\mu}=+1$ and $sign\frac{\partial\omega_{12}^*(\mu)}{\partial\mu}=+1$. Hence, case 1 of Table \ref{tab10} is proved. Treating similarly subcase (ii) of Case I and subcases (i) and (ii) of Case II, we can prove the rest of the cases in Table \ref{tab10}.\qed

{ {
The convergence conditions for the GMPSD(3) are given by the following.
\begin{cor}\label{cor:conv1b}
Consider the GMPSD(3) method. Let $A\in\mathbb{R}^{m\times m}$ and $Q\in\mathbb{R}^{n\times n}$ be
symmetric positive definite and $B\in\mathbb{R}^{m\times n}$ be of
full column rank.
Denote the minimum and the maximum eigenvalues of the matrix $J=Q^{-1}B^TA^{-1}B$ by $\mu_{min}$ and $\mu_{max}$, respectively. Then,
${\rho({\cg}{(\tau,\omega_1,\omega_2)})<1}$ if
\begin{equation}\label{eq:57}
    0<\tau<2,\ \ \omega_2<\omega_2^*(\mu_{max}) \ \ \mbox{and}\ \
    \omega_{13}^*(\mu_{max})<\omega_1<\omega_{14}^*(\mu_{max})
\end{equation}
where
\begin{equation}\label{eq:58}
\begin{array}{cccc}
\omega_2^*(\mu)=1-\frac{\tau^2\mu}{4},\;\;& \omega_{13}^*(\mu)=\frac{\tau-\omega_2}{1-\omega_2}-\frac{1}{\mu},\;\;&\omega_{14}^*(\mu)=\frac{2-\tau}{\tau\mu}+\frac{\tau-2\omega_2}{2(1-\omega_2)}.
\end{array}
\end{equation}
\end{cor}
{\bf{Proof}} Letting $a=0$ in the functional relationship (\ref{eq:26}) and following a similar approach as in the proof of Theorem \ref{theo:conv1bb}, we can prove (\ref{eq:57}).\qed

}}
Note that analogous results hold when $Q\in\mathbb{R}^{n\times n}$ is symmetric negative definite.

\subsection{Optimum parameters}
In the following theorem the optimum parameters of the GMPSD method are determined assuming that the matrix $Q$ is symmetric positive definite and $a\neq 0$.

\begin{theo}\label{theo:opt_1b}
Consider the GMPSD method. Let $A\in\mathbb{R}^{m\times m}$ and $Q\in\mathbb{R}^{n\times n}$ be symmetric positive definite and $B\in\mathbb{R}^{m\times n}$ be of
full column rank. Denote the minimum and the maximum eigenvalues of the matrix $J=Q^{-1}B^TA^{-1}B$ by
$\mu_{min}$ and $\mu_{max}$, respectively.
Then the spectral radius of the GMPSD method, $\rho(\cg(\tau_1, \tau_2, \omega_1, \omega_2,a))$, is minimized for any $\omega_2\neq\frac{\tau_{2_{opt}}}{\tau_{1_{opt}}}$ at
\begin{equation}\label{eq:optgmpsd11}
\omega_{1_{opt}}=\frac{\tau_{1_{opt}}(\tau_{2_{opt}}-\omega_2)}{\tau_{2_{opt}}-\tau_{1_{opt}}\omega_2}, \end{equation}
\begin{equation}\label{eq:optgmpsd12}
\tau_{1_{opt}}=\frac{4\sqrt{\mu_{min}\mu_{max}}}{(\sqrt{\mu_{min}}+\sqrt{\mu_{max}})^2}\ \ \mbox{and} \ \ \tau_{2_{opt}}=\frac{(1-a\omega_2)[1-(1-a)\omega_2]}{\sqrt{\mu_{min}\mu_{max}}}
\end{equation}
and its corresponding value is
\begin{equation}\label{eq:optgmpsd13}
\rho(\cg{(\tau_{1_{opt}},\tau_{2_{opt}},\omega_{1_{opt}},\omega_2, a)})=\frac{\sqrt{\mu_{max}}-\sqrt{\mu_{min}}}{\sqrt{\mu_{max}}+\sqrt{\mu_{min}}}.
\end{equation}
\end{theo}
{\bf{Proof}} Following a similar approach as in Theorem \ref{theo:opt_1}, using the functional relationship (\ref{eq:24}) and requiring $|\tilde\lambda_1|=|\tilde\lambda_N|$ we find
\begin{equation}\label{eq:optgmpsd134}
\tau_1\omega_1\omega_2-\tau_2\omega_1-\tau_1\omega_2+\tau_1\tau_2=0.
\end{equation}
Therefore, (\ref{eq:24}) because of (\ref{eq:optgmpsd134}), becomes
\begin{equation}\label{eq:optgmpsd135c}
\lambda^2+\lambda\left(\tau_1-2+\tau_1\hat\tau_2\mu\right)+1-\tau_1=0
\end{equation}
with
\begin{equation}\label{eq:optgmpsd135}
\hat\tau_2=\frac{\tau_2}{(1-a\omega_2)[1-(1-a)\omega_2]}
\end{equation}
which is the functional relationship of the GSOR method {{\cite{BaiParlWang2005}}}
with the only difference that now we have $(1-a\omega_2)[1-(1-a)\omega_2]$ instead of $1-a\omega_2$ in the denominator of $\hat\tau_2$ {{(see {(\ref{eq:5_s06bb5x1})})}}, hence (\ref{eq:optgmpsd11}), follows from (\ref{eq:optgmpsd134}) whereas (\ref{eq:optgmpsd12}) and (\ref{eq:optgmpsd13}) hold because of {{(\ref{eq:optgmpsd135c}), (\ref{eq:optgmpsd135}) and Theorem 4.1 in \cite{BaiParlWang2005}}}. \qed

\begin{cor}\label{cor:opt_1b}
Consider the simplified GMPSD method. Let $A\in\mathbb{R}^{m\times m}$ and $Q\in\mathbb{R}^{n\times n}$ be symmetric positive definite and $B\in\mathbb{R}^{m\times n}$ be of
full column rank. Denote the minimum and the maximum eigenvalues of the matrix $J=Q^{-1}B^TA^{-1}B$ by
$\mu_{min}$ and $\mu_{max}$, respectively. Then the spectral radius of the simplified GMPSD method, $\rho(\cg(\tau_1, \tau_2, \omega_1, 0, 0))$, is minimized at
\begin{equation}\label{eq:optgmpsd012b}
\omega_{1_{opt}}=\tau_{1_{opt}},
\end{equation}
\begin{equation}\label{eq:optgmpsd012c}
\tau_{1_{opt}}=\frac{4\sqrt{\mu_{min}\mu_{max}}}{(\sqrt{\mu_{min}}+\sqrt{\mu_{max}})^2}\ \ \mbox{and} \ \ \tau_{2_{opt}}=\frac{1}{\sqrt{\mu_{min}\mu_{max}}}
\end{equation}
and its corresponding value is
\begin{equation}\label{eq:optgmpsd013}
\rho(\cg{(\tau_{1_{opt}},\tau_{2_{opt}},\omega_{1_{opt}},0,0)})=\frac{\sqrt{\mu_{max}}-\sqrt{\mu_{min}}}{\sqrt{\mu_{max}}+\sqrt{\mu_{min}}}.
\end{equation}
\end{cor}
{\bf{Proof}} Letting $\omega_2=0$, (\ref{eq:optgmpsd11}), (\ref{eq:optgmpsd12}) and (\ref{eq:optgmpsd13}) yield (\ref{eq:optgmpsd012b}), (\ref{eq:optgmpsd012c}) and (\ref{eq:optgmpsd013}), respectively. \qed
It is worth noting here that the optimum values of $ \tau_{1_{opt}} ~ \mbox{and} ~ \tau_{2_{opt}} $ of the simplified GMPSD method are identical to the optimum values of $ \omega_{1_{opt}} ~ \mbox{and} ~ \omega_{2_{opt}} $ of the GSOR method, respectively.

{ {

\begin{cor}\label{cor:opt_1b3}
Consider the GMPSD(3) method. Let $A\in\mathbb{R}^{m\times m}$ and $Q\in\mathbb{R}^{n\times n}$ be symmetric positive definite and $B\in\mathbb{R}^{m\times n}$ be of
full column rank. Denote the minimum and the maximum eigenvalues of the matrix $J=Q^{-1}B^TA^{-1}B$ by
$\mu_{min}$ and $\mu_{max}$, respectively. If $\mmax<\frac{1}{4}$ or if $\mmax>\frac{1}{4}$ and either (i) $\mmin<\mu^*$ or (ii) $\mmin\geq\mu^*$ and $a_1\leq a\leq a_2$, then
the spectral radius of the GMPSD(3) method, $\rho(\cg{(\tau,\omega_{1},\omega_{2},a)})$, is minimized at
\begin{equation}\label{eq:optgmpsd31}
\omega_{1_{opt}}=\frac{\tau_{opt}-\omega_{2_{opt}}}{1-\omega_{2_{opt}}},
\end{equation}
\begin{equation}\label{eq:optgmssor32ca}
\tau_{{opt}}=\frac{4\sqrt{\mu_{min}\mu_{max}}}{(\sqrt{\mu_{min}}+\sqrt{\mu_{max}})^2} \;\;\mbox{and}\;\; \omega_{2_{opt}}=\frac{\sigma}{2\left[{1\pm\sqrt{1-a(1-a)\sigma}}\right]}
\end{equation}
and its corresponding value is
\begin{equation}\label{eq:optgmpsd33}
\rho(\cg{(\tau_{opt},\omega_{1_{opt}},\omega_{2_{opt}},a)})=\frac{\sqrt{\mu_{max}}-\sqrt{\mu_{min}}}{\sqrt{\mu_{max}}+\sqrt{\mu_{min}}}.
\end{equation}
where
\begin{equation}\label{eq:optgmpsd33bb}
\mu^*=\frac{\mmax}{(1-2\sqrt{\mmax})^2},\;\; a_1=\frac{2}{\sigma{{-}}\sqrt{\sigma(\sigma-4)}},\;\; a_2=\frac{2}{\sigma{{+}}\sqrt{\sigma(\sigma-4)}}
\end{equation}
with
\begin{equation}\label{eq:optgmpsd33bbc}
\sigma=4(1-M) \; \mbox{and}\; M=\displaystyle\frac{4\mmin\mmax}{(\sqrt{\mmin}+\sqrt{\mmax})^2}.
\end{equation}
\end{cor}
{\bf{Proof}}
Recall that GMPSD(3) is obtained by setting $\tau_1=\tau_2$ in GMPSD. Therefore, (\ref{eq:optgmpsd31}), (\ref{eq:optgmssor32ca}) and (\ref{eq:optgmpsd33}) are obtained by (\ref{eq:optgmpsd11}), (\ref{eq:optgmpsd12}) and (\ref{eq:optgmpsd13}) respectively.
In particular, by letting $\tau_{1_{opt}}=\tau_{2_{opt}}$ it follows from (\ref{eq:optgmpsd12}) that
\begin{equation}\label{eq:eq_alpha}
a(1-a)\omega_2^2-\omega_2+1-M=0,
\end{equation}
where $M$ is given by (\ref{eq:optgmpsd33bbc}).
This quadratic has real roots when
\begin{equation}\label{eq:D_thal}
a^2\sigma-a\sigma+1\geq0,
\end{equation}
where $\sigma$ is given by (\ref{eq:optgmpsd33bbc}).
Considering (\ref{eq:D_thal}) as a quadratic we distinguish two cases. Case 1: $\Delta_a<0$, Case 2: $\Delta_a\geq0$ where $\Delta_a=\sigma(\sigma-4)$.\\ \noindent
Case 1: $\Delta_a<0$. In this case we require $\sigma>0$ since $\sigma-4<0$ or in view of (\ref{eq:optgmpsd33bbc})
\begin{equation}\label{eq:optg3a2}
\sqrt{\mmin}(1-2\sqrt{\mmax})>-\sqrt{\mmax}.
\end{equation}
But, (\ref{eq:optg3a2}) holds if either $\mmax<\frac14$ or if $\mmax>\frac14$ and
$\mmin<\mu^*$ and (i) is proved. \\ \noindent
Case 2: $\Delta_a\geq0$. In this case we require $\sigma\leq0$ since $\sigma-4<0$ or,
because of (\ref{eq:optgmpsd33bbc}),
\begin{equation}\label{eq:optg3ab}
\sqrt{\mmin}(1-2\sqrt{\mmax})\leq -\sqrt{\mmax}
\end{equation}
which holds if $\mmax>\frac14$ and
\begin{equation}\label{eq:optg3abcb}
\mmin\geq\mu^*.
\end{equation}
In this case, for (\ref{eq:D_thal}) to hold, $a$ must lie in the range given by (ii).
Hence, the proof of the theorem is complete.\qed

\begin{cor}\label{cor:opt_1b4}
Consider the GMSSOR method. Let $A\in\mathbb{R}^{m\times m}$ and $Q\in\mathbb{R}^{n\times n}$ be symmetric positive definite and $B\in\mathbb{R}^{m\times n}$ be of
full column rank. Denote the minimum and the maximum eigenvalues of the matrix $J=Q^{-1}B^TA^{-1}B$ by
$\mu_{min}$ and $\mu_{max}$, respectively. If $\mmax<\frac{1}{4}$ or if $\mmax>\frac{1}{4}$ and either (i) $\mmin<\mu^*$ or (ii) $\mmin\geq\mu^*$ and $a_1\leq a\leq a_2$, then
the spectral radius of the
GMSSOR method, $\rho(\cg{(\omega_1, \omega_2, a)})$, is minimized at
\begin{equation}\label{eq:optgmssor32c}
\omega_{1_{opt}}=\frac{\hat\omega_{opt}-\omega_{2_{opt}}}{1-\omega_{2_{opt}}}
\end{equation}
where
\begin{equation}\label{eq:optgmssor31}
\hat\omega_{{opt}}=\frac{4\sqrt{\mu_{min}\mu_{max}}}{(\sqrt{\mu_{min}}+\sqrt{\mu_{max}})^2}\;\;\mbox{and}\;\;
\omega_{2_{opt}}=\frac{{{\sigma}}}{2\left[1\pm\sqrt{1-a(1-a)\sigma}\right]}
\end{equation}
and its corresponding value is
\begin{equation}\label{eq:optgmssor33}
\rho(\cg{(\omega_{1_{opt}},\omega_{2_{opt}},a)})=\frac{\sqrt{\mu_{max}}-\sqrt{\mu_{min}}}{\sqrt{\mu_{max}}+\sqrt{\mu_{min}}}.
\end{equation}
{{where $\mu^*, a_1, a_2, \sigma$ are given by (\ref{eq:optgmpsd33bb}), (\ref{eq:optgmpsd33bbc}).
}}
\end{cor}
{\bf{Proof}}
Recall that GMSSOR is obtained by setting $\tau_1=\tau_2=\hat\omega$ in GMPSD. Therefore, (\ref{eq:optgmssor32c}), (\ref{eq:optgmssor31}) and (\ref{eq:optgmssor33}) are obtained by (\ref{eq:optgmpsd11}), (\ref{eq:optgmpsd12}) and (\ref{eq:optgmpsd13}), respectively. Indeed, as in GMPSD(3), since $\tau_{1_{opt}}=\tau_{2_{opt}}$ it follows that (\ref{eq:eq_alpha})
holds also and by the analysis of the proof of Corollary \ref{cor:opt_1b3}, we have that (\ref{eq:optgmssor32c}), (\ref{eq:optgmssor31}) and (\ref{eq:optgmssor33})
hold under the same conditions as in Corollary \ref{cor:opt_1b3}.\qed

\begin{cor}\label{cor:opt_1b5}
{ {Consider the GPSD method. Let $A\in\mathbb{R}^{m\times m}$ and $Q\in\mathbb{R}^{n\times n}$ be symmetric positive definite and $B\in\mathbb{R}^{m\times n}$ be of
full column rank. Denote the minimum and the maximum eigenvalues of the matrix $J=Q^{-1}B^TA^{-1}B$ by
$\mu_{min}$ and $\mu_{max}$, respectively. If $\mmax<\frac{1}{4}$ or if $\mmax>\frac{1}{4}$ and either (i) $\mmin<\mu^*$ or (ii) $\mmin\geq\mu^*$ and $a_1\leq a\leq a_2$, then
the spectral radius of the GPSD  method, $\rho(\cg{(\tau, \omega, a)})$, is minimized at
\begin{equation}\label{eq:optgpsd32}
\tau_{{opt}}=\frac{4\sqrt{\mu_{min}\mu_{max}}}{(\sqrt{\mu_{min}}+\sqrt{\mu_{max}})^2} \;\; \mbox{and}\;\; \omega_{{opt}}=\frac{\sigma}{2\left[{1\pm\sqrt{1-a(1-a)\sigma}}\right]}
\end{equation}
and its corresponding value is
\begin{equation}\label{eq:optgpsd33bc0}
\rho(\cg{(\tau_{{opt}},\omega_{{opt}},a)})=\frac{\sqrt{\mu_{max}}-\sqrt{\mu_{min}}}{\sqrt{\mu_{max}}+\sqrt{\mu_{min}}}
\end{equation}
where
\begin{equation}\label{eq:optgmpsd33bb}
\mu^*=\frac{\mmax}{(1-2\sqrt{\mmax})^2},\;\; a_1=\frac{2}{\sigma{{-}}\sqrt{\sigma(\sigma-4)}},\;\; a_2=\frac{2}{\sigma{{+}}\sqrt{\sigma(\sigma-4)}}
\end{equation}
with
\begin{equation}\label{eq:optgmpsd33bbc}
\sigma=4(1-M) \; \mbox{and}\; M=\displaystyle\frac{4\mmin\mmax}{(\sqrt{\mmin}+\sqrt{\mmax})^2}.
\end{equation}
}}
\end{cor}
{\bf{Proof}}{ {
GPSD follows from GMPSD by letting $\tau=\tau_1=\tau_2$ and $\omega=\omega_1=\omega_2$ or $\tau_{opt}=\tau_{1_{opt}}=\tau_{2_{opt}}$ and $\omega_{opt}=\omega_{1_{opt}}=\omega_{2_{opt}}$. By equating the expressions of $\tau_{1_{opt}}$ and $\tau_{2_{opt}}$ given by (\ref{eq:optgmpsd12}) we obtain
\begin{equation}\label{eq:eq_alpha}
a(1-a)\omega_2^2-\omega_2+1-M=0,
\end{equation}
where $M$ is given by (\ref{eq:optgmpsd33bbc}).
This quadratic has real roots when
\begin{equation}\label{eq:D_thal}
a^2\sigma-a\sigma+1\geq0,
\end{equation}
where $\sigma$ is given by (\ref{eq:optgmpsd33bbc}).
We distinguish two cases. Case 1: $\Delta_a<0$, Case 2: $\Delta_a\geq0$ where $\Delta_a=\sigma(\sigma-4)$.\\ \noindent
Case 1: $\Delta_a<0$. In this case we require $\sigma>0$ since $\sigma-4<0$ or in view of (\ref{eq:optgmpsd33bbc})
\begin{equation}\label{eq:optg3a2}
\sqrt{\mmin}(1-2\sqrt{\mmax})>-\sqrt{\mmax}.
\end{equation}
But, (\ref{eq:optg3a2}) holds if either $\mmax<\frac14$ or if $\mmax>\frac14$ and
$\mmin<\mu^*$ hence (i) is proved. \\ \noindent
Case 2: $\Delta_a\geq0$. In this case we require $\sigma\leq0$ since $\sigma-4<0$ or,
because of (\ref{eq:optgmpsd33bbc}),
\begin{equation}\label{eq:optg3ab}
\sqrt{\mmin}(1-2\sqrt{\mmax})\leq -\sqrt{\mmax}
\end{equation}
which holds if $\mmax>\frac14$ and
\begin{equation}\label{eq:optg3abcb}
\mmin\geq\mu^*.
\end{equation}
In this case, for (\ref{eq:D_thal}) to hold, $a$ must lie in the range given by (ii).
Therefore, it follows that $\omega_{opt}$ is given by (\ref{eq:optgpsd32}).
\qed
}}

}}
Analogous results hold in case where the matrix Q is symmetric negative definite.

\section{Numerical results}
In this section we study the numerical solution of the following linear Stokes equation
\begin{equation}\label{eq:stokes}
\begin{cases}
-\mu\Delta{\bf{u}}+ \nabla w=\tilde f,& \mbox{στο}\; \Omega \\
\hspace{1.2cm}\nabla\cdot{\bf{u}}=\tilde g,& \mbox{στο}\; \Omega\\
\hspace{1.75cm}{\bf{u}}=0, & \mbox{στο}\; \partial\Omega\\
\hspace{0.4cm}\int_{\Omega}w(x)dx=0,&
\end{cases}
\end{equation}
{{where $\Omega=(0,1)\times(0,1)\subset\mathbb{R}^2$, $\partial\Omega$ is the boundary of $\Omega$, $\Delta$ is the componentwise Laplace operator, {\bf{u}} is a vector-valued function representing the velocity and $w$ is a scalar function representing the pressure. Furthermore, we assume that the functions $\tilde f, \tilde g$ are constant. By discretizing  (\ref{eq:stokes}) with the upwind scheme, we obtain the system of linear equations (\ref{eq:01}), in which  \cite{BaiGolPan2004}
\begin{equation*}
A=\left(
\begin{array}{cc}
 I\otimes T+T\otimes I & 0\\
0 & I\otimes T+T\otimes I
\end{array}
\right) \in \mathbb{R}^{2p^2\times 2p^2},
\end{equation*}
\begin{equation*}
B=\left(
\begin{array}{c}
 I\otimes F\\
F\otimes I
\end{array}
\right)\in \mathbb{R}^{2p^2\times p^2}
\end{equation*}
\noindent with \noindent
\begin{equation*}
T=\frac{\mu}{h^2}\cdot tridiag(-1,2,-1)\in \mathbb{R}^{p\times p}, \;\;\; F=\frac{1}{h}\cdot tridiag(-1,1,0)\in \mathbb{R}^{p\times p},
\end{equation*}
\noindent  $h=\frac{1}{p+1}$ being the discretization mesh size and $\otimes$ the Kronecker product symbol.}}
For this example, we let $\mu=1$, $m=2p^2$ and $n=p^2$.
Hence, the total number of variables is $m+n=3p^2$. \\ \noindent
We choose the matrix $Q$ to be an approximation to $B^TA^{-1}B$. The reason being that if $Q\simeq B^TA^{-1}B$ then $J=Q^{-1}B^TA^{-1}B\simeq I$. In this case the ratio of the maximum to the minimum eigenvalue of the matrix $J$ becomes minimum and its value is approximately 1. As a consequence, the spectral radius of the iteration matrix of the GMESOR and GMPSD methods attains its minimum value. We choose $Q$, according to the following two cases:\\\\ \noindent
1.  $Q=B^T\hat A^{-1}B,\; \hat A=\mbox{tridiag}(A)$\\ \noindent
2.  $Q=B^T\hat A^{-1}B,\; \hat A=\mbox{diag}(A)$,\\\\ \noindent
\noindent where $\hat A$ is the tridiagonal or the diagonal part of $A$. The choice of the matrix $\hat A$ instead of $A$ is due to the difficulty in computing the inverse matrix of $A$.
In this example the eigenvalues of $Q$ are real and positive.
{{In actual computations, we choose the right-hand-side vector $(b^T, q^T)^T\in\mathbb{R}^{m+n}$ such that the exact solution of the augmented linear system (\ref{eq:01}) is
$({(x^{*})}^T, {(y^{*})}^T)^T=(1,1,...,1)^T\in\mathbb{R}^{m+n}$, and perform all runs in MATLAB (version $R{\_}2012b$) with a machine precision $10^{-16}$. The machine used was an Intel i5 personal computer with 6G memory.
In our computations, all runs are started from the initial vector $\left({(x^{(0)})}^T, {y^{(0)}}^T\right)^T=0$, and terminated if the current iterations satisfy
$$RES=\frac{\sqrt{||b-Ax^{(k)}-By^{(k)}||^2_2+||q-B^Tx^{(k)}||^2_2}}{\sqrt{||b-Ax^{(0)}-By^{(0)}||^2_2+||q-B^Tx^{(0)}||^2_2}}\leq 10^{-9},$$ where $RES$ is the norm of absolute residual vectors, or if the numbers of the prescribed iterations $k_{max}=1200$ are exceeded. We also use the same example to compare our methods with the PHSS  \cite{BaiGolPan2004} {{and Krylov subspace methods \cite{SaSch86}, \cite{Saad03}, \cite{HVdV03}}}.
}}
\\ \noindent In Table \ref{tab3} we computed the optimal parameters $\tau_{1_{opt}}, \tau_{2_{opt}}\; \mbox{and}\; \omega_{2_{opt}}$ and the optimal spectral radius $\rho_{opt}$ of the GMESOR method, for various problem sizes (m,n) using (\ref{eq:optgmesor1}), (\ref{eq:optgmesor2}) and (\ref{eq:optgmesor3}).
Furthermore, we computed the  optimum parameters {{$\tau_{2_{opt}}({exp})$, $\omega_{2_{opt}}({exp})$ }} and the spectral radii  {{$\rho(\tau_{2_{opt}}{(exp)})$ and $\rho(\omega_{2_{opt}}{(exp)})$}},  experimentally by trial and error. The parameter $\tau_{1} $ was kept fixed and was given its optimum value.
Our results show that {{$\rho_{opt}\simeq\rho(\tau_{2_{opt}}{(exp)})\simeq\rho(\omega_{2_{opt}}{(exp)})$}}
and {{$\omega_{2_{opt}}=\tau_{2_{opt}}\simeq\tau_{2_{opt}}({exp})\simeq\omega_{2_{opt}}({exp})$}} thus verifying  Theorem \ref{theo:opt_1}.
\noindent {{The numerical results in Table \ref{tab5} verify that the parameter $a$ may be chosen arbitrary, while the minimum value of $\rho(\ch(\tau_{1}, \tau_2, \omega_{2},a))$ remains approximately the same.  $\rho(\ch(\tau_{1_{opt}}, \tau_{2_{opt}}, \omega_{2_{opt}},a))$ was computed using Matlab. The slightly different values are due to rounding errors.}}
{{Finally, in Table {\ref{tab6}} we list numerical results with respect to the number of total iteration steps (denoted by ``ITER"), the elapsed CPU time in seconds (denoted by ``CPU") and $RES$ for the GSOR, GMESOR and Simplified GMPSD iterative methods.
We remark that our numerical results verify the validity of theorem {\ref{theo:opt_1}} and corollary {\ref{cor:opt_1b}}, since GSOR, GMESOR and Simplified GMPSD methods require the same number of iterations for convergence. Indeed, this was expected since all these methods have the same spectral radius for the optimum values of their parameters. Note that all the aforementioned methods require approximately the same computing time.
Furthermore, for comparison purposes we also considered the {{PHSS($a^*$)}}, GMRES, GMRES($\#$), PGMRES and PGMRES($\#$) methods. The integer $\#$ in GMRES($\#$) and PGMRES($\#$) methods denotes the number of restarting steps, while the integer $a^*$ denotes the theoretical optimal parameter of the PHSS method. We also list numerical results with respect to the number of total iteration steps and the elapsed CPU time in seconds for these methods. The preconditioned matrix $Q$ in {{PHSS($a^*$)}} is given by the aforementioned cases 1 and 2. The preconditioner, say $K$, for the PGMRES and PGMRES($\#$) methods is given by \cite{ElmSch86}, \cite{ElmSilWa02}, \cite{MurGoWa00}, \cite{WaSil93}
\begin{equation*}
K=\left[
\begin{array}{cc}
 \hat A & 0\\
0 & I
\end{array}
\right].
\end{equation*}
We remark that the GSOR, GMESOR and Simplified GMPSD methods always outperform the other testing methods, except of the {{PHSS($a^*$)}} method, considerably with respect to iteration steps as $p$ increases. However, the overall computing time of the GSOR, GMESOR and Simplified GMPSD methods is much smaller than that of all the other testing methods.
With * we denote that the method converges but after too many hours.  With regard to the matrix $Q$, Case 1 is the best choice for all methods tested as it requires the least iteration steps and CPU times.
}}
\begin{table}[htbp]
\caption{Experimental results showing that $\omega_{2_{opt}}=\tau_{2_{opt}}$ in GMESOR. } \label{tab3}
\begin{center}
\begin{tabular}{|c|c|c|c|c|c|c} \hline
\multicolumn{2}{|c|}{\bf{m}}&\multicolumn{1}{|c|}{128}&\multicolumn{1}{|c|}{512}&\multicolumn{1}{|c|}{1152}\\
\hline\multicolumn{2}{|c|}{\bf{n}}&\multicolumn{1}{|c|}{64}&\multicolumn{1}{|c|}{256}&\multicolumn{1}{|c|}{576}\\
\hline\multicolumn{2}{|c|}{\bf{m+n}}&\multicolumn{1}{|c|}{192}&\multicolumn{1}{|c|}{768}&\multicolumn{1}{|c|}{1728}\\
\hline
& $\tau_{1_{opt}}$&0.663309&0.442911&0.330674\\
\cline{2-5}    &$\tau_{2_{opt}}$& 0.499375&0.285422&0.198468\\
\cline{2-5}   &$\omega_{2_{opt}}$& 0.499375&0.285422&0.198468\\
\cline{2-5}{\bf {Case 1}}   &$\rho_{{opt}}$& 0.580251&0.746384&0.811229\\
\cline{2-5}     &{{$\tau_{2_{opt}}(exp)$}}&0.5 &0.286&0.199\\
\cline{2-5}    &{{$\rho(\tau_{2_{opt}}{(exp)})$}}&0.582936 &0.750508&0.823517\\
\cline{2-5}    &{{$\omega_{2_{opt}}({exp})$}}& 0.499&0.285&0.198\\
\cline{2-5}   &{{$\rho(\omega_{2_{opt}}{(exp)})$}}& 0.581866&0.749401&0.822877\\
\hline
& $\tau_{1_{opt}}$&0.757767 &0.631420&0.558518\\
\cline{2-5}    &$\tau_{2_{opt}}$&1.950825&2.529944&2.974309\\
\cline{2-5}   &$\omega_{2_{opt}}$& 1.950825&2.529944&2.974309\\
\cline{2-5}{\bf {Case 2}}   &$\rho_{{opt}}$&  0.492171&0.607108&0.664441\\
\cline{2-5}     &{{$\tau_{2_{opt}}(exp)$}}&1.951&2.530&2.975\\
\cline{2-5}    &{{$\rho(\tau_{2_{opt}}{(exp)})$}}&0.492374&0.607155&0.664925 \\
\cline{2-5}    &{{$\omega_{2_{opt}}({exp})$}}& 1.950&2.529&2.974\\
\cline{2-5}   &{{$\rho(\omega_{2_{opt}}{(exp)})$}}&0.493127 &0.607901&0.664657\\
\hline
\end{tabular}
\end{center}
\end{table}
\begin{table}[htbp]
\caption{\footnotesize{{Computation of $\rho(\ch(\tau_{1_{opt}}, \tau_{2_{opt}}, \omega_{2_{opt}},a))$ for various values of the parameter $a$ (Case 1, $p=40$). }}} \label{tab5}
\begin{center}{{
\begin{tabular}{|c|c|c|c|c|c|c|c|c|c|c} \hline
\cline{1-4}    $a$  &  $\tau_{1_{opt}}$     &   $\tau_{2_{opt}}=\omega_{2_{opt}}$    &   $\rho(\ch(\tau_{1_{opt}}, \tau_{2_{opt}}, \omega_{2_{opt}},a))$  \\
\cline{1-4}     0   &   2.18851E-001    &   1.229935E-001   &   0.883807    \\ \hline
\cline{1-4}     10  &   2.18851E-001    &   5.515564E-002   &   0.883808    \\ \hline
\cline{1-4}    $10^2$ &   2.18851E-001    &   9.248083E-003   &   0.883808     \\ \hline
\cline{1-4}    $10^3$ &   2.18851E-001    &   9.919351E-004   &   0.883809   \\ \hline
\cline{1-4}    $10^4$ &   2.18851E-001    &   9.991876E-005   &   0.883810     \\ \hline
\cline{1-4}    $10^5$ &   2.18851E-001    &   9.999187E-006   &   0.883807   \\ \hline
\cline{1-4}    $10^6$ &   2.18851E-001    &   9.999919E-007   &   0.883807       \\ \hline
\cline{1-4}    $10^7$ &   2.18851E-001    &   9.999992E-008   &   0.883807   \\ \hline
\cline{1-4}    $10^8$ &   2.18851E-001    &   9.999999E-009   &   0.883815   \\ \hline
\cline{1-4}    $10^9$ &   2.18851E-001    &   1.000000E-009   &   0.886425    \\ \hline
\cline{1-4}  $10^{10}$&   2.18851E-001    &   1.000000E-010   &   0.883879     \\ \hline
\end{tabular}
}}\end{center}
\end{table}

\begin{table}[htbp]
{{\caption{ITER, CPU and RES for the testing methods} \label{tab6}}}
\begin{center}
{{\scriptsize{\begin{tabular}{|c|c|c|c|c|c|c|c|c|c|c||c|c|c|c|} \hline
\cline{3-9}\multicolumn{2}{|c|}{}&\multicolumn{1}{|c|}{\bf{p}}&\multicolumn{1}{|c|}{8}&\multicolumn{1}{|c|}{16}&\multicolumn{1}{|c|}{24}&\multicolumn{1}{|c|}{32}&\multicolumn{1}{|c|}{40}&\multicolumn{1}{|c|}{48}\\
\cline{3-9}\multicolumn{2}{|c|}{}&\multicolumn{1}{|c|}{\bf{n}}&\multicolumn{1}{|c|}{64}&\multicolumn{1}{|c|}{256}&\multicolumn{1}{|c|}{576}&\multicolumn{1}{|c|}{1024}&\multicolumn{1}{|c|}{1600}&\multicolumn{1}{|c|}{2304}\\
\cline{3-9}\multicolumn{2}{|c|}{}&\multicolumn{1}{|c|}{\bf{m}}&\multicolumn{1}{|c|}{128}&\multicolumn{1}{|c|}{512}&\multicolumn{1}{|c|}{1152}&\multicolumn{1}{|c|}{2048}&\multicolumn{1}{|c|}{3200}&\multicolumn{1}{|c|}{4608}\\
\cline{3-9}\multicolumn{2}{|c|}{}&\multicolumn{1}{|c|}{\bf{m+n}}&\multicolumn{1}{|c|}{192}&\multicolumn{1}{|c|}{768}&\multicolumn{1}{|c|}{1728}&\multicolumn{1}{|c|}{3072}&\multicolumn{1}{|c|}{4800}&\multicolumn{1}{|c|}{6912}\\
\hline\hline
&& ITER & 46 & 86 & 126 & 167 & 207 & 248\\
\cline{3-9}&GSOR& CPU & 0,05 & 0,36 & 3,71 & 22,47 & 86,47 & 258,28\\
\cline{3-9}&& RES & 6,79E-10 & 9,04E-10 & 9,79E-10 & 8,97E-10 & 9,74E-10 & 9,44E-10\\
\cline{2-9}&& ITER & 46 & 86 & 126 & 167 & 207 & 248\\
\cline{3-9}&GMESOR& CPU & 0,05 & 0,36& 3,71 & 22,61 & 86,78 & 258,93\\
\cline{3-9}&& RES & 6,79E-10 & 9,04E-10 & 9,79E-10 & 8,97E-10 & 9,74E-10 & 9,44E-10\\
\cline{2-9}& & ITER & 46 & 86 & 126 & 167 & 207 & 248\\
\cline{3-9}& Simplified & CPU & 0,05 & 0,35& 3,71 & 22,59 & 86,54 & 258,28\\
\cline{3-9}&GMPSD & RES & 7,03E-10 & 9,12E-10 & 9,83E-10 & 8,99E-10 & 9,75E-10 & 9,45E-10\\
\cline{2-9}& & ITER & 24 & 35 & 44 & 51 & 57 & 63\\
\cline{3-9}{Case 1}&PHSS$(a^*)$& CPU & 0,34 & 5,24& 34,52 & 147,95 & 472,87 & 1247,80\\
\cline{3-9}&& RES & 6,19E-10 & 9,62E-10 & 7,63E-10 & 7,36E-10 & 9,63E-10 & 8,82E-10\\
\cline{2-9}& & ITER & 73 & 176 & 285 & 386 & 506 & 606\\
\cline{3-9}& GMRES& CPU & 0,33 & 9,24& 155,67 & 1.240,62 & 6.352,04 & 22.214,42\\
\cline{2-9}& & ITER & 73 & 327 & 831 & 1794 & 3436 & 9965\\
\cline{3-9}&GMRES(100)& CPU & 0,26 & 16,12 & 404,70 & 5.417,99 & 41.626,66 & 356.831,22\\
\cline{2-9}& & ITER & 76 & 143 & 207 & 275 & 344 & 410\\
\cline{3-9}& PGMRES & CPU & 0,50 & 11,35 & 130,19 & 956,58 & 4.557,29 & 15.684,49\\
\cline{2-9}&& ITER & 76 & 178 & 321 & 509 & 1038 & 1281\\
\cline{3-9}& PGMRES(100) & CPU & 0,36 & 11,02& 172,80 & 1.615,68 & 12.838,90 & 46.595,50\\
\hline\hline
&& ITER & 65 & 124 & 182 & 241 & 300 & 359\\
\cline{3-9}&GSOR& CPU & 0,07 & 0,42 & 4,11 & 24,55 & 93,24 & 278,67\\
\cline{3-9}&& RES & 8,35E-10 & 8,25E-10 & 9,32E-10 & 9,14E-10 & 9,19E-10 & 9,35E-10\\
\cline{2-9}&& ITER & 65 & 124 & 182 & 241 & 300 & 359\\
\cline{3-9}&GMESOR& CPU & 0,06 & 0,40 & 4,08 & 24,41 & 93,14 & 278,48\\
\cline{3-9}&& RES & 8,35E-10 & 8,25E-10 & 9,32E-10 & 9,14E-10 & 9,19E-10 & 9,35E-10\\
\cline{2-9}&& ITER & 65 & 124 & 182 & 241 & 300 & 359\\
\cline{3-9}&Simplified& CPU & 0,07 & 0,39 & 4,11 & 24,51 & 93,20 & 278,08\\
\cline{3-9}&GMPSD& RES & 8,55E-10 & 8,30E-10 & 9,35E-10 & 9,15E-10 & 9,20E-10 & 9,35E-10\\
\cline{2-9}&& ITER & 29 & 43 & 53 & 62 & 69 & 76\\
\cline{3-9}{Case 2}&PHSS$(a^*)$& CPU & 0,34 & 5,28 & 34,64 & 148,51 & 474,44 & 1261,96\\
\cline{3-9}&& RES & 9,88E-10 & 6,53E-10 & 7,99E-10 & 8,48E-10 & 9,61E-10 & 9,72E-10\\
\cline{2-9}& & ITER & 73 & 176 & 285 & 386 & 506 & 606\\
\cline{3-9}& GMRES& CPU & 0,33 & 9,24& 155,67 & 1.240,62 & 6.352,04 & 22.214,42\\
\cline{2-9}& & ITER & 73 & 327 & 831 & 1794 & 3436 & 9965\\
\cline{3-9}&GMRES(100)& CPU & 0,26 & 16,12 & 404,70 & 5.417,99 & 41.626,66 & 356.831,22\\
\cline{2-9}& & ITER & 75 & 164 & 253 & 347 & 446 & 537\\
\cline{3-9}& PGMRES & CPU & 0,50 & 12,18 & 151,43 & 1.168,31 & 5.761,78 & 20.131,99\\
\cline{2-9}&& ITER & 75 & 275 & 544 & 997 & * & *\\
\cline{3-9}& PGMRES(100)& CPU & 0,36 & 15,31 & 279,17 & 3.072,70 & $\gg$ 72h & $\gg$ 72h \\
\hline
\end{tabular}}
}}\end{center}
\end{table}

\section{Remarks and Conclusions}
In this paper we studied the impact of two different preconditioning matrices on the convergence of iterative methods for the solution of the augmented linear system (\ref{eq:01}) when the coefficient matrix $\ca$ is of the form (\ref{eq:02}). We assumed that $A\in \mathbb{R}^{m\times m}$ was a symmetric positive definite matrix and $B\in\mathbb{R}^{m\times n}$ was a matrix of full column rank, where $m\geq n$, whereas $Q$ was a symmetric positive or negative definite matrix.
Under these assumptions we were able to find
sufficient conditions for the GMESOR and GMPSD iterative methods to converge. Further, using a geometric analysis analogous to Varga \cite{Var1} we determined the optimum values of the parameters of all methods studied such as to attain the maximum rate of convergence.  From our analysis it was shown that GMESOR and GMPSD are equivalent since they have the same spectral radius for the optimum values of their parameters, which is given by {{ (\ref{eq:optgmesor3}).}}
This result was verified by our numerical experiments, where the simplified GMPSD, the GMESOR and the GSOR methods require approximately the same computing time. Moreover, all the aforementioned methods outperform the {{PHSS($a^*$)}}, {{GMRES, GMRES($\#$), PGMRES and PGMRES($\#$) methods}} considerably with respect to CPU times. It is worth mentioning that, for the saddle point problem, the GMPSD method has a similar behavior as the Modified PSD (MPSD) method for two-cyclic matrices \cite{LoukMis2}. Indeed, in \cite{LoukMis2} we proved the equivalence of MPSD and MSOR methods for two-cyclic matrices in case the eigenvalues of the Jacobi matrix are either all real or all imaginary. However, it is believed that this equivalence will not hold for the case where the eigenvalues of the $J$ matrix are complex.

\begin{acknowledgements}
The authors would like to thank the referees for their constructive comments and suggestions which improved considerably the original form of the paper. The second author would like to thank the Department of Applied Mathematics and Statistics of State University of New York at Stony Brook for its warm hospitality while working on the paper.
\end{acknowledgements}



\end{document}